\documentclass[a4paper,10pt]{article}
\usepackage[reqno]{amsmath}
\makeatletter
\newcommand{\leqnomode}{\tagsleft@true}
\newcommand{\reqnomode}{\tagsleft@false}
\makeatother

\usepackage[english,british]{babel}
\usepackage{amsmath,amsthm,amssymb,dsfont,CJK,esint}
\usepackage{adjustbox,lipsum}
\usepackage{tikz}
\usepackage{standalone}
\usepackage{mathtools}
\usepackage{fancyhdr}
\usepackage{graphicx}
\usepackage{bibentry}
\usepackage{enumitem}
\usepackage{xcolor}
\usepackage{stackengine}
\usepackage[
            pdfstartview=FitH,
            CJKbookmarks=true,
            bookmarksnumbered=true,
            bookmarksopen=true,
            colorlinks,
            pdfborder=001,
            linkcolor=red,
            citecolor=green,
            ]{hyperref}
\usepackage[a4paper,width=150mm,top=30mm,bottom=30mm]{geometry}
\hypersetup{
    colorlinks,
    citecolor=green,
    filecolor=black,
    linkcolor=red,
    urlcolor=black
}

\newtheorem{theorem}{Theorem}[section]
\newtheorem{lemma}[theorem]{Lemma}
\newtheorem{definition}[theorem]{Definiton}
\newtheorem{proposition}[theorem]{Proposition}

\theoremstyle{definition}
\newtheorem{remx}[theorem]{Remark}
\newenvironment{remark}
  {\pushQED{\qed}\remx}
  {\popQED\endremx}

\makeatletter
\newsavebox\myboxA
\newsavebox\myboxB
\newlength\mylenA

\newcommand*\yoverline[2][0.75]{%
    \sbox{\myboxA}{$\m@th#2$}%
    \setbox\myboxB\null
    \ht\myboxB=\ht\myboxA%
    \dp\myboxB=\dp\myboxA%
    \wd\myboxB=#1\wd\myboxA
    \sbox\myboxB{$\m@th\overline{\copy\myboxB}$}
    \setlength\mylenA{\the\wd\myboxA}
    \addtolength\mylenA{-\the\wd\myboxB}%
    \ifdim\wd\myboxB<\wd\myboxA%
       \rlap{\hskip 0.5\mylenA\usebox\myboxB}{\usebox\myboxA}%
    \else
        \hskip -0.5\mylenA\rlap{\usebox\myboxA}{\hskip 0.5\mylenA\usebox\myboxB}%
    \fi}
\makeatother

\numberwithin{equation}{section}

\begin{document}




\newcommand{\diver}{\operatorname{div}}
\newcommand{\lin}{\operatorname{Lin}}
\newcommand{\curl}{\operatorname{curl}}
\newcommand{\ran}{\operatorname{Ran}}
\newcommand{\kernel}{\operatorname{Ker}}
\newcommand{\la}{\langle}
\newcommand{\ra}{\rangle}
\newcommand{\N}{\mathbb{N}}
\newcommand{\R}{\mathbb{R}}
\newcommand{\C}{\mathbb{C}}

\newcommand{\ld}{\lambda}
\newcommand{\fai}{\varphi}
\newcommand{\0}{0}
\newcommand{\n}{\mathbf{n}}
\newcommand{\uu}{{\boldsymbol{\mathrm{u}}}}
\newcommand{\UU}{{\boldsymbol{\mathrm{U}}}}
\newcommand{\buu}{\bar{{\boldsymbol{\mathrm{u}}}}}
\newcommand{\ten}{\\[4pt]}
\newcommand{\six}{\\[4pt]}
\newcommand{\nb}{\nonumber}
\newcommand{\hgamma}{H_{\Gamma}^1(\OO)}
\newcommand{\opert}{O_{\varepsilon,h}}
\newcommand{\barx}{\bar{x}}
\newcommand{\barf}{\bar{f}}
\newcommand{\hatf}{\hat{f}}
\newcommand{\xoneeps}{x_1^{\varepsilon}}
\newcommand{\xh}{x_h}
\newcommand{\scaled}{\nabla_{1,h}}
\newcommand{\scaledb}{\widehat{\nabla}_{1,\gamma}}
\newcommand{\vare}{\varepsilon}
\newcommand{\A}{{\bf{A}}}
\newcommand{\RR}{{\bf{R}}}
\newcommand{\B}{{\bf{B}}}
\newcommand{\CC}{{\bf{C}}}
\newcommand{\D}{{\bf{D}}}
\newcommand{\K}{{\bf{K}}}
\newcommand{\oo}{{\bf{o}}}
\newcommand{\id}{{\bf{Id}}}
\newcommand{\E}{\mathcal{E}}
\newcommand{\ii}{\mathcal{I}}
\newcommand{\sym}{\mathrm{sym}}
\newcommand{\lt}{\left}
\newcommand{\rt}{\right}
\newcommand{\ro}{{\bf{r}}}
\newcommand{\so}{{\bf{s}}}
\newcommand{\e}{{\bf{e}}}
\newcommand{\ww}{{\boldsymbol{\mathrm{w}}}}
\newcommand{\vv}{{\boldsymbol{\mathrm{v}}}}
\newcommand{\zz}{{\boldsymbol{\mathrm{z}}}}
\newcommand{\U}{{\boldsymbol{\mathrm{U}}}}
\newcommand{\G}{{\boldsymbol{\mathrm{G}}}}
\newcommand{\VV}{{\boldsymbol{\mathrm{V}}}}
\newcommand{\T}{{\boldsymbol{\mathrm{U}}}}
\newcommand{\II}{{\boldsymbol{\mathrm{I}}}}
\newcommand{\ZZ}{{\boldsymbol{\mathrm{Z}}}}
\newcommand{\hKK}{{{\bf{K}}}}
\newcommand{\f}{{\bf{f}}}
\newcommand{\g}{{\bf{g}}}
\newcommand{\lkk}{{\bf{k}}}
\newcommand{\tkk}{{\tilde{\bf{k}}}}
\newcommand{\W}{{\boldsymbol{\mathrm{W}}}}
\newcommand{\Y}{{\boldsymbol{\mathrm{Y}}}}
\newcommand{\EE}{{\boldsymbol{\mathrm{E}}}}
\newcommand{\F}{{\bf{F}}}
\newcommand{\spacev}{\mathcal{V}}
\newcommand{\spacevg}{\mathcal{V}^{\gamma}(\Omega\times S)}
\newcommand{\spacevb}{\bar{\mathcal{V}}^{\gamma}(\Omega\times S)}
\newcommand{\spaces}{\mathcal{S}}
\newcommand{\spacesg}{\mathcal{S}^{\gamma}(\Omega\times S)}
\newcommand{\spacesb}{\bar{\mathcal{S}}^{\gamma}(\Omega\times S)}
\newcommand{\skews}{H^1_{\barx,\mathrm{skew}}}
\newcommand{\kk}{\mathcal{K}}
\newcommand{\OO}{O}
\newcommand{\bhe}{{\bf{B}}_{\vare,h}}
\newcommand{\pp}{{\mathbb{P}}}
\newcommand{\ff}{{\mathcal{F}}}
\newcommand{\mWk}{{\mathcal{W}}^{k,2}(\Omega)}
\newcommand{\mWa}{{\mathcal{W}}^{1,2}(\Omega)}
\newcommand{\mWb}{{\mathcal{W}}^{2,2}(\Omega)}
\newcommand{\twos}{\xrightharpoonup{2}}
\newcommand{\twoss}{\xrightarrow{2}}
\newcommand{\bw}{\bar{w}}
\newcommand{\br}{\bar{{\bf{r}}}}
\newcommand{\bz}{\bar{{\bf{z}}}}
\newcommand{\tw}{{W}}
\newcommand{\tr}{{{\bf{R}}}}
\newcommand{\tz}{{{\bf{Z}}}}
\newcommand{\lo}{{{\bf{o}}}}
\newcommand{\hoo}{H^1_{00}(0,L)}
\newcommand{\ho}{H^1_{0}(0,L)}
\newcommand{\hotwo}{H^1_{0}(0,L;\R^2)}
\newcommand{\hooo}{H^1_{00}(0,L;\R^2)}
\newcommand{\hhooo}{H^1_{00}(0,1;\R^2)}
\newcommand{\dsp}{d_{S}^{\bot}(\barx)}
\newcommand{\LB}{{\bf{\Lambda}}}
\newcommand{\LL}{\mathbb{L}}
\newcommand{\mL}{\mathcal{L}}
\newcommand{\mhL}{\widehat{\mathcal{L}}}
\newcommand{\loc}{\mathrm{loc}}
\newcommand{\tqq}{\mathcal{Q}^{*}}
\newcommand{\tii}{\mathcal{I}^{*}}
\newcommand{\Mts}{\mathbb{M}}
\newcommand{\pot}{\mathrm{pot}}
\newcommand{\tU}{{\widehat{\bf{U}}}}
\newcommand{\tVV}{{\widehat{\bf{V}}}}
\newcommand{\pt}{\partial}
\newcommand{\bg}{\Big}
\newcommand{\hA}{\widehat{{\bf{A}}}}
\newcommand{\hB}{\widehat{{\bf{B}}}}
\newcommand{\hCC}{\widehat{{\bf{C}}}}
\newcommand{\hD}{\widehat{{\bf{D}}}}
\newcommand{\fder}{\partial^{\mathrm{MD}}}
\newcommand{\Var}{\mathrm{Var}}
\newcommand{\pta}{\partial^{0\bot}}
\newcommand{\ptaj}{(\partial^{0\bot})^*}
\newcommand{\ptb}{\partial^{1\bot}}
\newcommand{\ptbj}{(\partial^{1\bot})^*}
\newcommand{\geg}{\Lambda_\vare}
\newcommand{\tpta}{\tilde{\partial}^{0\bot}}
\newcommand{\tptb}{\tilde{\partial}^{1\bot}}
\newcommand{\ua}{u_\alpha}
\newcommand{\pa}{p\alpha}
\newcommand{\qa}{q(1-\alpha)}
\newcommand{\Qa}{Q_\alpha}
\newcommand{\Qb}{Q_\eta}
\newcommand{\ga}{\gamma_\alpha}
\newcommand{\gb}{\gamma_\eta}
\newcommand{\ta}{\theta_\alpha}
\newcommand{\tb}{\theta_\eta}


\newcommand{\mH}{\mathcal{H}}
\newcommand{\mD}{\mathcal{D}}
\newcommand{\csob}{\mathcal{S}}
\newcommand{\gna}{{\mathrm{C}_{\mathrm{GN}}}}
\newcommand{\gnb}{{\widehat{\mathrm{C}}_{\mathrm{GN}}}}
\newcommand{\mA}{\mathcal{A}}
\newcommand{\mK}{\mathcal{K}}
\newcommand{\mS}{\mathcal{S}}
\newcommand{\mI}{\mathcal{I}}
\newcommand{\mM}{\mathcal{M}}
\newcommand{\tas}{{4}}
\newcommand{\tbs}{{6}}
\newcommand{\tass}{{0}}
\newcommand{\tbss}{{\fraco}}
\newcommand{\tm}{{\tilde{m}}}
\newcommand{\tdu}{{\tilde{u}}}
\newcommand{\tpsi}{{\tilde{\psi}}}
\newcommand{\Z}{{\mathds{Z}}}
\newcommand{\fraco}{{\frac{1}{2}}}
\newcommand{\fraca}{{\frac{1}{4}}}
\newcommand{\fracb}{{\frac{1}{6}}}
\newenvironment{proof1}{\paragraph{Proof of Theorem \ref{main theorem 1}}}{\hfill$\square$}

\title{Sharp scattering threshold for the cubic-quintic NLS in the focusing-focusing regime}
\author{Yongming Luo \thanks{Institut f\"{u}r Wissenschaftliches Rechnen, Technische Universit\"at Dresden, Germany} \thanks{\href{mailto:yongming.luo@tu-dresden.de}{Email: yongming.luo@tu-dresden.de}}
}

\date{}
\maketitle

\begin{abstract}
We consider the large data scattering problem for the 2D and 3D cubic-quintic nonlinear Schr\"odinger equation in the focusing-focusing regime. Our attention is firstly restricted to the 2D space, where the cubic nonlinearity is $L^2$-critical. We establish a new type of scattering criterion that is uniquely determined by the mass of the initial data, which differs from the classical setting based on the Lyapunov functional. At the end, we formulate a solely mass-determining scattering threshold for the 3D cubic-quintic nonlinear Schr\"odinger equation in the focusing-focusing regime.
\end{abstract}

\section{Introduction and main results}
\label{intro}
In this paper, we consider the cubic-quintic nonlinear Schr\"odinger equation
\begin{equation}\label{NLS0}
i\pt_t u+\Delta u+\mu_1|u|^{2}u+\mu_2|u|^{4}u=0\quad\text{in $\R\times\R^d$}
\end{equation}
for $d=2,3$. The cubic-quintic nonlinear Schr\"odinger equation (CQNLS) serves as a toy model in many physical applications such as nonlinear optics and Bose-Einstein condensation. Physically, the cubic and quintic nonlinearities model the two-body and three-body interactions respectively. The signs $\mu_i$ can be tuned to be defocusing ($\mu_i<0$) or focusing ($\mu_i>0$), indicating the repulsivity or attractivity of the many-body interactions. We refer to \cite{phy2,phy3,phy1} and the references therein for a comprehensive introduction on the physical background of the CQNLS.

On the other hand, the CQNLS has also attracted much attention from the mathematical community due to its abundant analytical structure: one easily verifies the $L^2$-criticality of the quintic term in 1D, the $L^2$-criticality of the cubic term in 2D and the $\dot{H}^1$-criticality of the quintic term in 3D. Additionally, the mixed type nature of the CQNLS prevents any possible application of scaling invariance property, which makes the mathematical analysis more subtle and challenging. For recent mathematical progress on the study of the CQNLS with a particular focus on the scattering and blow-up phenomenon, we refer to \cite{TaoVisanZhang,MiaoDoubleCrit,killip_visan_soliton,Cheng2020,Carles_Sparber_2021,Murphy2021CPDE,killip2020cubicquintic,carles2020soliton}.

Our particular interest is firstly devoted to establishing a sharp scattering threshold for the 2D CQNLS in the focusing-focusing regime ($\mu_1,\mu_2>0$), which has not been considered in any of the above mentioned references. In particular, we see no possibility to easily adapt the existing arguments to the focusing-focusing model; To achieve our aim, we need some new ideas and ingredients. For simplicity, we set $\mu_1=\mu_2=1$ and in the following we consider the normalized CQNLS
\begin{equation}\label{NLS}
i\pt_t u+\Delta u+|u|^{2}u+|u|^{4}u=0\quad\text{in $\R\times\R^2$}.
\end{equation}
In \cite{TaoVisanZhang}, it was shown that as long as $\mu_1,\mu_2<0$ (namely both nonlinearities are defocusing), \eqref{NLS0} is globally well-posed and scatters in time for any initial data $u(0)=u_0\in H^1(\R^2)$. The proof relies on the so-called interaction Morawetz inequalities and the defocusing nature of the model is essential, hence the proof can not be adapted to other types of models. In fact, the result from \cite{TaoVisanZhang} does not hold when at least one of the $\mu_i$ is positive: \eqref{NLS0} might possess solutions that blow-up in finite time, or soliton solutions. Here, the soliton solutions are referred to solutions $u$ of \eqref{NLS0} having the form $u(t,x)=e^{i\omega t}S(x)$ with $\omega\in\R$, where $S$ satisfies the stationary CQNLS
\begin{align}
-\Delta S+\omega S-\mu_1|S|^{2}S-\mu_2|S|^{4}S=0\quad\text{in $\R^2$}.\label{standing wave eq0}
\end{align}
For the normalized focusing-focusing model, \eqref{standing wave eq0} reads
\begin{align}\label{standing wave eq}
-\Delta S+\omega S-|S|^{2}S-|S|^{4}S=0\quad\text{in $\R^2$}.
\end{align}
It turns out that the energy level corresponding to soliton solutions is the minimal threshold for a solution being failed to meet the dichotomy of scattering and blow-up. Our starting point is the result in \cite{SoaveSubcritical} given by Soave, where the author studied the existence problem of solutions of \eqref{standing wave eq}. Therein, Soave considered the following variational problem
\begin{align}\label{soave problem}
m_c=\inf_{u\in H^1(\R^2)}\{\mH(u):\mM(u)=c,\mK(u)=0\}
\end{align}
for $c>0$, where
\begin{align*}
\mM(u)&:=\|u\|^2_2,\nonumber\\
\mH(u)&:=\fraco\|\nabla u\|_2^2-\fraca\|u\|^{\tas}_{\tas}-\fracb\|u\|^{\tbs}_{\tbs},\nonumber\\
\mK(u)&:=\|\nabla u\|_2^2-\fraco\|u\|_\tas^\tas-\frac{2}{3}\|u\|_\tbs^\tbs.
\end{align*}
Physically, $\mM(u),\mH(u),\mK(u)$ denote the mass, energy and virial respectively. It was shown in \cite{SoaveSubcritical} that $\mK(u)=0$ is a natural constraint, thus using the Lagrange multiplier theorem we know that any optimizer of \eqref{soave problem} is automatically a solution of \eqref{standing wave eq}. An optimizer of $m_c$ is also said to be a \textit{ground state} since it has the least energy among all candidates. To formulate the result in \cite{SoaveSubcritical}, we also denote by $Q$ the unique positive and radially symmetric solution of
\begin{align*}
-\Delta Q+Q-Q^{3}=0.
\end{align*}
Having all the preliminaries we are able to introduce the following result from \cite{SoaveSubcritical}:
\begin{theorem}[\cite{SoaveSubcritical}]\label{soave}
We have the following existence and blow-up results:
\begin{itemize}
\item[(i)]\textbf{Existence of ground state}: For any $c\in(0,{\mM(Q)})$ the variational problem \eqref{soave problem} has a minimizer $S_c$ with $\mH(S_c)=m_c\in(0,\infty)$. Moreover, $S_c$ is a solution of \eqref{standing wave eq} with some $\omega>0$. In addition, $S_c$ can be chosen to be positive and radially symmetric.

\item[(ii)]\textbf{Blow-up criterion}: Assume that $u_0\in H^1(\R^2)$ satisfies the conditions $\mM(u_0)\in(0,{\mM(Q)})$, $\mH(u_0)<m_{\mM(u_0)}$ and $\mK(u_0)<0$. Assume also that $|x|u_0\in L^2(\R^2)$. Then the solution $u$ of \eqref{NLS} with $u(0)=u_0$ blows-up in finite time.
\end{itemize}
\end{theorem}
The aim of the present paper is to show that the threshold given by Soave is exactly the sharp scattering threshold for \eqref{NLS}.
\begin{theorem}\label{main theorem 1}
Define the set
\begin{align}
\mA&:=\{u\in H^1(\R^2):\mM(u)<{\mM(Q)},\mH(u)<m_{\mM(u)},\mK(u)> 0\}\label{scattering threshold}
\end{align}
and assume that $u_0\in \mA$. Then the solution $u$ of \eqref{NLS} with $u(0)=u_0$ is global and scatters in time.
\end{theorem}


We should compare the scattering criterion \eqref{scattering threshold} with the one given in \cite{Akahori2013}, which is nowadays the golden rule for large data scattering problems of NLS with combined power type nonlinearities. Therein, the authors considered the NLS
\begin{align}\label{aka nls}
i\pt_t u+\Delta u+|u|^{p-2} u+|u|^{\frac{4}{d-2}}u=0\quad\text{in $\R\times\R^d$}
\end{align}
for $d\geq 5$ and $p\in(2+\frac{4}{d},2+\frac{4}{d-2})$, namely the focusing energy-critical NLS with a focusing mass-supercritical and energy-subcritical perturbation. To formulate the scattering threshold, the authors imposed the so-called Lyapunov functional
\begin{align}
\mS_\omega(u):=\frac{\omega}{2} \mM(u)^2+\mH(u)
\end{align}
and considered the variational problem
\begin{align}\label{gamma omega}
\gamma_\omega:=\inf_{u\in H^1(\R^2)\setminus\{0\}}\{\mS_\omega(u):\mK(u)=0\}.
\end{align}
The following result is due to \cite{Akahori_soliton,Akahori2013,akahori2019global}:
\begin{theorem}[\cite{Akahori_soliton,Akahori2013,akahori2019global}]
Let $d\geq 3$ and $\omega>0$. Then
\begin{itemize}
\item[(i)] For any $d\geq 3$ and $\omega>0$ we have $\gamma_\omega\in(0,\frac{\csob^{\frac{d}{2}}}{d}]$, where $\csob$ is the optimal constant for the Sobolev inequality.
\item[(ii)]For any $d\geq 4$ and $\omega>0$ and for $d=3$ and any sufficiently small $\omega>0$ we have $\gamma_\omega\in(0,\frac{\csob^{\frac{d}{2}}}{d})$. Moreover, the variational problem \eqref{gamma omega} possesses an optimizer $P_\omega$. Consequently, the optimizer $P_\omega$ is a soliton solution of \eqref{aka nls} with the given $\omega$.
\item[(iii)] Assume that
\begin{align}
u_0\in\{v\in H^1(\R^d):\mS_\omega(v)<\gamma_\omega,\mK(v)> 0\}.\label{aka scatter criterion}
\end{align}
Additionally we assume that $u_0$ is radially symmetric when $d=3$. Then the solution $u$ of \eqref{aka nls} with $u(0)=u_0$ is global and scatters in time.
\end{itemize}
\end{theorem}
We should point out that the scattering results given in \cite{Akahori2013} were originally formulated in the cases $d\geq 5$. The results can nonetheless be extended to all dimensions $d\geq 3$ in a natural way by combining with the results from \cite{KenigMerle2006} and the recent published paper \cite{Dodson4dfocusing}. The scattering criterion \eqref{aka scatter criterion} has been later successfully applied in \cite{Miao3dradial,Xie2014,Miao_5d_2016,Miao_4d,MiaoDoubleCrit,Xu_Yang_2016} to formulate a sharp scattering threshold for NLS with combined power type nonlinearities in different regimes. However, \eqref{aka scatter criterion} seems not to be compatible with problems having focusing $L^2$-critical nonlinearity due to the following reason: The constraint $\mM(u)<\mM(Q)$ is essential since the $L^2$-critical nonlinearity $|u|^{\frac{4}{d}}u$ shares the same scaling of the Laplacian $\Delta u$ under the scaling operator
$$ T_\ld u(x):=\ld^{\frac{d}{2}}u(\ld x).$$
However, by considering the Lyapunov functional $\mS_\omega$ we would have double constraints on the mass of the initial data, which might violate the conciseness of the scattering threshold. By comparison with \eqref{aka scatter criterion} we also see that the scattering criterion \eqref{scattering threshold} has the advantage that it is uniquely determined by the mass of the initial data. Such a solely mass-determining threshold could also be more physically relevant in the following sense: besides being a conserved quantity, the mass also measures many physically important quantities such as the power supply in nonlinear optics, or the total number of particles in the Bose-Einstein condensation.

As a surprising byproduct, we are able to formulate a scattering criterion for the problem \eqref{aka nls} within the framework of the present paper by also invoking the results from \cite{SoaveCritical,wei2021normalized}. In particular, there is no smallness condition and mass constraint in all dimensions and the energy threshold is exactly the ground state energy which is positive and smaller than $d^{-1}\csob^{\frac{d}{2}}$. We will continue the discussion in Section \ref{outlook section}. As a consequence of Theorem \ref{3d scattering threshold} given below, we are able to impose the following solely mass-determining scattering threshold for the 3D CQNLS in the focusing-focusing regime:
\begin{theorem}\label{3d scattering thm}
Let $d=3$ and $\mu_1=\mu_2=1$. Define the set
\begin{align}
\mathcal{B}&:=\{u\in H^1(\R^3):\mH(u)<m_{\mM(u)},\mK(u)> 0\}
\end{align}
(where $\mK(u)$ is suitably redefined in the 3D case) and assume that $u_0\in \mathcal{B}\cap H^1_{\mathrm{rad}}(\R^3)$. Then the solution $u$ of \eqref{NLS0} with $u(0)=u_0$ is global and scatters in time.
\end{theorem}

\begin{remark}
We note that unlike the 2D case, in Theorem \ref{3d scattering thm} there is no mass constraint imposed for the initial data. This difference stems from the fact that the cubic nonlinearity is mass-critical and mass-supercritical in 2D and 3D respectively. To be more precise, when applying the $L^2$-scaling $u\mapsto T_\ld u$ to the quantity $\frac12 \|\nabla u\|_2^2-\frac14 \|u\|_4^4$ in 2D, we see that
\begin{align*}
\frac12 \|\nabla (T_\ld u)\|_2^2-\frac14 \|T_\ld u\|_4^4=\ld^2\bg(\frac12 \|\nabla u\|_2^2-\frac14 \|u\|_4^4\bg).
\end{align*}
Therefore, the quantity $\frac12 \|\nabla u\|_2^2-\frac14 \|u\|_4^4$ does not vary w.r.t. scaling and keeping the mass below the (mass-critical) ground state is essential for applications of Gagliardo-Nirenberg inequalities. Such heuristics do not hold any longer in the 3D case. Indeed, by applying the $L^2$-scaling to the quantity $\frac12 \|\nabla u\|_2^2-\frac14 \|u\|_4^4$ in 3D we obtain
\begin{align*}
\frac12 \|\nabla (T_\ld u)\|_2^2-\frac14 \|T_\ld u\|_4^4=\frac{\ld^2}{2} \|\nabla u\|_2^2-\frac{\ld^3}{4} \|u\|_4^4.
\end{align*}
We see in this case that the quantity $\frac12 \|\nabla u\|_2^2-\frac14 \|u\|_4^4$ does not play the same role as in 2D and the mass constraint is no longer relevant. Rather, we should consider the function $\ld\mapsto \frac{\ld^2}{2} \|\nabla u\|_2^2-\frac{\ld^3}{4} \|u\|_4^4$ and the corresponding variational analysis becomes more delicate in comparison with the 2D case. We refer to the papers \cite{SoaveSubcritical,SoaveCritical,wei2021normalized} for a more detailed survey on such phenomenon.
\end{remark}

The study on the existence and stability results of \eqref{standing wave eq0} is also a very interesting topic. In this direction, we refer to the classical papers \cite{lions1,Berestycki1983,weinstein,Cazenave_orbital,jeanjean_classical,brezis_nirenberg} and also the recent papers \cite{BellazziniJeanjean2016,SoaveSubcritical,SoaveCritical,wei2021normalized}.

\subsubsection*{Roadmap for the proof of Theorem \ref{main theorem 1}}
We summarize here briefly the idea for the proof of Theorem \ref{main theorem 1}. The proof follows the classical concentration compactness arguments given by Kenig and Merle \cite{KenigMerle2006}: by assuming that the claim in Theorem \ref{main theorem 1} does not hold, we are able to derive a minimal blow-up solution $u_c$ of \eqref{NLS} with
$$ \|\la\nabla\ra^{\frac{1}{2}}u_c\|_{L_{t,x}^4(\R)}=\infty,$$
which also satisfies $u_c=0$ and thus leads to a contradiction. The main challenge arises from the fact that since the scattering is considered w.r.t. the $H^1$-topology, it is impossible to prove Theorem \ref{main theorem 1} only relying on the energy $\mH$. It is at this point to note that the $H^1$-norm of a solution $u$ of \eqref{NLS} can be controlled by the Lyapunov functional $\mS_\omega(u)$, which is not the case here. To build up the inductive contradiction hypothesis, we should rather take the mass and energy of the initial data into account simultaneously. To be more precise, we
utilize the mass-energy-indicator functional (MEI-functional) $\mD$ introduced in \cite{killip_visan_soliton} to derive a minimal blow-up solution. The idea can be described as follows: a mass-energy pair $(\mM(u),\mH(u))$ being admissible implies $\mD(u)\in(0,\infty)$; In order to escape the admissible region $\Omega$, a function $u$ must approach the boundary of $\Omega$ and one would deduce $\mD(u)\to\infty$. We can therefore assume that the supremum $\mD^*$ of $\mD(u)$ running over all admissible $u$ is finite, which leads to a contradiction and we conclude that $\mD^*=\infty$, which will finish the desired proof. However, the situation by the focusing-focusing model is more delicate: a mass-energy pair being admissible does not automatically imply the positivity of the virial $\mK$. In particular, it is not trivial that the linear profiles would have positive virial at the first glance. We will appeal to the geometric properties of the MEI-functional $\mD$, combining with the variational arguments from \cite{Akahori2013}, to overcome this difficulty.

\subsubsection*{Outline of the paper}
In Section \ref{Auxiliary preliminaries} we collect some auxiliary tools from \cite{Cheng2020,Murphy2021CPDE} that will be useful by the construction of the minimal blow-up solution; Section \ref{Variational problems and estimates} is devoted to the variational estimates and the construction of the MEI-functional $\mD$; Finally, we prove in Section \ref{Existence of the minimal blow-up solution} and Section \ref{Extinction of the minimal blow-up solution} the existence and extinction of the minimal blow-up solution respectively. In Section \ref{outlook section} we formulate a scattering criterion for the problem \eqref{aka nls} under the framework of the present paper. In \ref{m0 mQ} we establish the precise endpoint values $m_0$ and $m_Q$ of the curve $c\mapsto m_c$ in 2D.

\subsection{Notation and definitions}\label{Notations and definitions}
We use the notation $A\lesssim B$ whenever there exists some positive constant $C$ such that $A\leq CB$. Similarly we define $A\gtrsim B$ and  use $A\sim B$ when $A\lesssim B\lesssim A$. We denote by $\|\cdot\|_p$ the $L^p(\R^d)$-norm for $p\in[1,\infty]$. We similarly define the $H^1(\R^d)$-norm by $\|\cdot\|_{H^1}$. The following quantities will be used throughout the paper:
\begin{align*}
\mM(u)&:=\|u\|^2_2,\nonumber\\
\mH(u)&:=\fraco\|\nabla u\|_2^2-\fraca\|u\|^{\tas}_{\tas}-\fracb\|u\|^{\tbs}_{\tbs},\nonumber\\
\mK(u)&:=\|\nabla u\|_2^2-\fraco\|u\|_\tas^\tas-\frac{2}{3}\|u\|_\tbs^\tbs,\nonumber\\
\mI(u)&:=\mH(u)-\frac{1}{2}\mK(u)=\frac{1}{6}\|u\|_\tbs^\tbs.
\end{align*}
We will also frequently use the scaling operator
\begin{align*}
T_\ld u(x)&:=\ld^{\frac{d}2} u(\ld x).
\end{align*}
One easily verifies that the $L^2$-norm is invariant under this scaling. We denote by $Q$ the unique positive and radially symmetric ground state of
\begin{align*}
-\Delta Q+Q-Q^{3}=0.
\end{align*}
For the existence and uniqueness of $Q$, we refer to \cite{weinstein} and \cite{Kwong_uniqueness} respectively. We denote by $\mathrm{C}_{\mathrm{GN}}$ the 2D optimal $L^2$-critical Gagliardo-Nirenberg constant, i.e.
\begin{align}\label{def gn l2 crit}
\mathrm{C}_{\mathrm{GN}}=\inf_{u\in H^1(\R^2)\setminus\{0\}}\frac{\|\nabla u\|_2^2\|u\|_2^{2}}{\|u\|_{\tas}^{\tas}}.
\end{align}
Using Pohozaev identities (see for instance \cite{lions1}) and scaling arguments one easily verifies that
\begin{align}\label{GN-L1}
\mathrm{C}_{\mathrm{GN}}=\frac{1}{2}\mM(Q).
\end{align}
We also denote by $\widehat{\mathrm{C}}_{\mathrm{GN}}$ the optimal Gagliardo-Nirenberg constant for the quintic nonlinearity, i.e.
\begin{align}
\widehat{\mathrm{C}}_{\mathrm{GN}}=\inf_{u\in H^1(\R^2)\setminus\{0\}}\frac{\|\nabla u\|_2^4\|u\|_2^{2}}{\|u\|_{\tbs}^{\tbs}}.
\end{align}
For $d\geq 3$ we denote by $\csob$ the optimal constant of the Sobolev inequality, i.e.
\begin{align*}
\csob:=\inf_{u\in\mathcal{D}^{1,2}(\R^d)\setminus \{0\}}\frac{\|\nabla u\|_2^2}{\|u\|_{2^*}^2}.
\end{align*}
Here, the space $\mathcal{D}^{1,2}(\R^d)$ is defined by
\begin{align*}
\mathcal{D}^{1,2}(\R^d):=\{u\in L^{2^*}(\R^d):\nabla u\in L^2(\R^d)\}
\end{align*}
and $2^*=\frac{2d}{d-2}$. For an interval $I\subset \R$, the space $L_t^qL_x^r(I)$ is defined by
\begin{align*}
L_t^qL_x^r(I):=\{u:I\times \R^2\to\C:\|u\|_{L_t^qL_x^r(I)}<\infty\},
\end{align*}
where
\begin{align*}
\|u\|^q_{L_t^qL_x^r(I)}:=\int_{\R}\|u\|^q_r\,dt.
\end{align*}
A pair $(q,r)$ is said to be $L^2$-admissible in 2D if $q,r\in[2,\infty]$, $\frac{2}{q}+\frac{2}{r}=1$ and $(q,r)\neq(2,\infty)$. For any $L^2$-admissible pairs $(q_1,r_1)$ and $(q_2,r_2)$ we have the following Strichartz estimates: if $u$ is a solution of
\begin{align}
i\pt_t u+\Delta u=F(u)
\end{align}
in $I\subset\R$ with $t_0\in I$ and $u(t_0)=u_0$, then
\begin{align}
\|u\|_{L_t^q L_x^r(I)}\lesssim \|u_0\|_2+\|F(u)\|_{L_t^{q_2'} L_x^{r_2'}(I)},
\end{align}
where $(q_2',r_2')$ is the H\"older conjugate of $(q_2,r_2)$. For a proof, we refer to \cite{EndpointStrichartz,Cazenave2003}. The $S$-norm is defined by
\begin{align}
\|u\|_{S(I)}:=\sup_{q\in(2^+,\infty]}\|u\|_{L_t^q L_x^r(I)},
\end{align}
where the supremum is taken over all $L^2$-admissible pairs $(q,r)$ with $q\in(2^+,\infty]$ and $2^+>2$ is some positive constant that is sufficiently close to $2$. In this paper, the scattering concept is referred to the following definition:
\begin{definition}[Scattering]
A global solution $u$ of \eqref{NLS} is said to be forward in time scattering if there exists some $\phi_+\in H^1(\R^2)$ such that
\begin{align}
\lim_{t\to\infty}\|u(t)-e^{it\Delta}\phi_+\|_{H^1}=0.
\end{align}
A backward in time scattering solution is similarly defined. $u$ is then called a scattering solution when it is both forward and backward in time scattering.
\end{definition}
We define the Fourier transform of a function $f$ by
\begin{align*}
\hat{f}(\xi)=\mathcal{F}(f)(\xi)=\int_{\R^2}f(x)e^{-i\xi\cdot x}\,dx.
\end{align*}
For $s\in\R$, the multipliers $|\nabla|^s$ and $\la\nabla\ra^s$ are defined by the symbols
\begin{align*}
|\nabla|^s f(x)&=\mathcal{F}^{-1}\bg(|\xi|^s\hat{f}(\xi)\bg)(x),\\
\la\nabla\ra^s f(x)&=\mathcal{F}^{-1}\bg((1+|\xi|^2)^{\frac{s}{2}}\hat{f}(\xi)\bg)(x).
\end{align*}
Let $\psi\in C^\infty_c(\R^2)$ be a fixed radial, non-negative function such that $\psi(x)=1$ if $|x|\leq 1$ and $\psi(x)=0$ for $|x|\geq \frac{11}{10}$. Then for $N>0$, we define the Littlewood-Paley projectors by
\begin{align*}
P_{\leq N} f(x)&=\mathcal{F}^{-1}\bg(\psi\bg(\frac{\xi}{N}\bg)\hat{f}(\xi)\bg)(x),\\
P_{N} f(x)&=\mathcal{F}^{-1}\bg(\bg(\psi\bg(\frac{\xi}{N}\bg)-\psi\bg(\frac{2\xi}{N}\bg)\bg)\hat{f}(\xi)\bg)(x),\\
P_{\geq N} f(x)&=\mathcal{F}^{-1}\bg(\bg(1-\psi\bg(\frac{\xi}{N}\bg)\bg)\hat{f}(\xi)\bg)(x).
\end{align*}

\section{Auxiliary preliminaries}\label{Auxiliary preliminaries}
In this section we collect some useful auxiliary lemmas from \cite{Cheng2020,Murphy2021CPDE}. The proofs will be omitted here and we refer to \cite{Cheng2020,Murphy2021CPDE} for further details. We begin with the small data well-posedness and stability theory, which can be proved in a standard way.

\begin{lemma}[Small data well-posedness]\label{well posedness lemma}
For any $A>0$ there exists some $\beta>0$ such that the following is true: Suppose that $I$ is some interval and $t_0\in I$. Suppose also that $u_0\in H^1(\R^2)$ with
\begin{align}
\|u_0\|_{H^1}\leq A
\end{align}
and
\begin{align}
\| \la\nabla\ra e^{i(t-t_0)\Delta}u_0\|_{{L_{t,x}^4}(I)}\leq \beta.
\end{align}
Then \eqref{NLS} has a unique solution $u\in C(I;H^1(\R^2))$ with $u(t_0)=u_0$ such that
\begin{align}
\|\la\nabla\ra u\|_{S(I)}&\lesssim \|u_0\|_{H^1},\\
\|\la\nabla\ra u\|_{{L_{t,x}^4}(I)}&\leq 2\|\la\nabla\ra e^{i(t-t_0)\Delta}u_0\|_{{L_{t,x}^4}(I)}.
\end{align}
Denote by $I_{\max}$ the maximal lifespan of $u$. We then have the following blow-up and scattering criterion: if the solution $u$ of \eqref{NLS} satisfies
\begin{align}
\|\la\nabla\ra^{\frac{1}{2}} u\|_{{L_{t,x}^4}(I_{\max})}<\infty,
\end{align}
then $I_{\max}=\R$ and $u$ scatters in both positive and negative time.
\end{lemma}

\begin{remark}\label{remark}
Using Strichartz we infer that
\begin{align}
\|\la\nabla\ra e^{i(t-t_0)\Delta}u_0\|_{{L_{t,x}^4}(I)}\lesssim \|u_0\|_{H^1}.
\end{align}
Thus Lemma \ref{well posedness lemma} is applicable for all $u_0$ with sufficiently small $H^1$-norm.
\end{remark}

\begin{lemma}[Stability]\label{long time pert}
Let $u\in C(I;H^1(\R^2))$ be a solution of \eqref{NLS} defined on some interval $I\ni t_0$. Assume also that $w$ is an approximate solution of the following perturbed NLS
\begin{align}
i\pt_t w+\Delta w=-|w|^{2}w-|w|^{4}w+e
\end{align}
such that
\begin{align}
\|w\|_{L_t^\infty H_x^1(I)}&\leq B_1,\label{condition c1}\\
\|\la\nabla\ra^{\frac{1}{2}}(u(t_0)-w(t_0))\|_{2}&\leq B_2,\label{condition c3}\\
\|\la\nabla\ra^{\frac{1}{2}}w\|_{{L_{t,x}^4}(I)}&\leq B_3\label{condition c2}
\end{align}
for some $B_1,B_2,B_3>0$. Then there exists some positive $\beta_0=\beta_0(B_1,B_2,B_3)\ll 1$ with the following property: if
\begin{align}
\|\la\nabla\ra^{\frac{1}{2}} e^{i(t-t_0)\Delta}(u(t_0)-w(t_0))\|_{{L_{t,x}^4}(I)}&\leq \beta,\label{condition a}\\
\|\la\nabla\ra^{\frac{1}{2}} e\|_{L_{t,x}^{\frac{4}{3}}(I)}&\leq\beta\label{condition b}
\end{align}
for some $0<\beta<\beta_0$, then
\begin{align}
\|\la\nabla \ra^{\frac{1}{2}} (u-w)\|_{{L_{t,x}^4}(I)}&\leq C(B_1,B_2,B_3)\beta,\\
\|\la\nabla \ra^{\frac{1}{2}} (u-w)\|_{S(I)}&\leq C(B_1,B_2,B_3),\\
\|\la\nabla \ra^{\frac{1}{2}} u\|_{S(I)}&\leq C(B_1,B_2,B_3).
\end{align}
\end{lemma}

Next we introduce the linear profile decomposition used in present paper. Since \eqref{NLS} is a focusing mass-critical NLS with a mass-supercritical and energy-subcritical perturbation, we should apply an $L^2$-profile decomposition on the approximating sequence $(\psi_n)_n$ rather than an $\dot{H}^1$-profile decomposition as in \cite{Keraani_h1}. The classical $L^2$-profile decomposition was originally given by \cite{l2_decomp_R,merle_vega_l2_decomp,l2_all} and later applied in \cite{MiaoDoubleCrit} for the radial mass-energy double critical NLS. To remove the radial restriction we should appeal to the following linear profile decomposition from \cite{Cheng2020}:

\begin{lemma}[Linear profile decomposition]\label{linear profile}
Let $(\psi_n)_n$ be a bounded sequence in $H^1(\R^2)$. Then up to a subsequence of $(\psi_n)_n$, there exist some number $K^*\in\N\cup\{\infty\}$, a sequence of nonzero linear profiles $(\phi^j)_{1\leq j\leq K^*}\subset L^2(\R^2)$, a sequence of symmetry parameters $(\ld_n^j,t_n^j,x_n^j,\xi_n^j)_{n\in\N,1\leq j\leq K^*}\subset(0,\infty)\times\R\times \R^2\times\R^2$ with $\sup_{n\in\N}|\xi_n^j|\lesssim_j 1$ and a sequence of remainders $(w_{n}^k)_{n\in\N,1\leq k\leq K^*}\subset H^1(\R^2)$ such that
\begin{itemize}
\item[(i)] The parameters $(\ld_n^j,t_n^j,x_n^j,\xi_n^j)_{n,j}$ satisfy
\begin{align}
\lim_{n\to\infty}&\bg\{\bg|\log \frac{\ld_n^j}{\ld_n^l}\bg|+\frac{|t_n^j-t_n^l|}{(\ld_n^j)^2}\nonumber\\
&\quad\quad\quad+\ld_n^j|\xi_n^j-\xi_n^l|+\frac{|x_n^j-x_n^l+2t_n^j(\xi_n^j-\xi_n^l)|}{\ld_n^j}\bg\}=\infty
\end{align}
for all finite $1\leq j,l\leq K^*$ with $j\neq l$. Moreover,
\begin{align}
\lim_{n\to\infty}\ld_n^j&=\ld^j_\infty\in\{1,\infty\},\\
\ld_n^j&\equiv1\quad\text{if $\ld_\infty^j=1$}.
\end{align}

\item[(ii)] There exists some $\theta\in(0,1)$ such that for any finite $1\leq k\leq K^*$ we have the decomposition
\begin{align}
\psi_n=\sum_{j=1}^k T_n^j P_n^j \phi^j+w_n^k,
\end{align}
where
\begin{align}
T_n^j u(x):=e^{ix\cdot\xi_n^j}e^{-it_n^j\Delta}(\ld_n^j)^{-1}u\bg((\ld_n^j)^{-1}(\cdot-x_n^j)\bg)(x)
\end{align}
and
\begin{equation}
P_n^j\phi^j=\left\{
             \begin{array}{ll}
             \phi^j,&\text{if $\ld_\infty^j=1$},\\
             P_{\leq(\ld_n^j)^\theta}\phi^j,&\text{if $\ld_\infty^j=\infty$}.
             \end{array}
\right.
\end{equation}
Moreover, if $\ld^j_\infty=1$, then $\xi_n^j\equiv 0$ and $\phi^j\in H^1(\R^2)$.

\item[(iii)] The remainders $(w_n^k)_{n,k}$ satisfy
\begin{align}
\lim_{k\to K^*}\lim_{n\to \infty}\|\la\nabla\ra^{\frac{1}{2}} e^{it\Delta}w_n^k\|_{{L_{t,x}^4}(\R)}=0.
\end{align}

\item[(iv)] The following orthogonal properties are satisfied: for $s\in\{0,1\}$ and finite $1\leq J\leq K^*$ we have
\begin{align}
\||\nabla|^s\psi_n\|_2^2&=\sum_{j=1}^J \||\nabla|^s T_n^j P_n^j\phi^j\|_2^2+\||\nabla |^s w_n^J\|_2^2+o_n(1),\\
\mH(\psi_n)&=\sum_{j=1}^J \mH(T_n^j P_n^j\phi^j)+\mH( w_n^J)+o_n(1),\\
\mI(\psi_n)&=\sum_{j=1}^J \mI(T_n^j P_n^j\phi^j)+\mI( w_n^J)+o_n(1),\\
\mK(\psi_n)&=\sum_{j=1}^J \mK(T_n^j P_n^j\phi^j)+\mK( w_n^J)+o_n(1).
\end{align}
\end{itemize}
\end{lemma}

The following large scale approximation result is an immediate consequence of Lemma \ref{long time pert} and the scattering result in \cite{Dodson4dmassfocusing} for focusing mass-critical NLS.

\begin{lemma}[Large scale approximation]\label{large scale proxy}
Let $\phi^j$ be a linear profile given through Lemma \ref{linear profile}. Suppose also that $\mM(\phi^j)<\mM(Q)$ and $\ld_\infty^j= \infty$. Then for all sufficiently large $n$ (possibly depending on $j$) there exists a solution $v^j_n$ of \eqref{NLS} such that $v_n^j$ is a global and scattering solution with $v^j_n(0)=T_n^jP_n^j\phi^j$. Moreover, we have
\begin{align}
\limsup_{n\to\infty}\|\la\nabla \ra^{\frac12}v_n^j\|_{L_{t,x}^4(\R)}\lesssim_{\mM(\phi^j)}1.
\end{align}
\end{lemma}

\section{Variational estimates}\label{Variational problems and estimates}
In this section we derive some variational estimates as preliminaries for the proofs given in Section \ref{Existence of the minimal blow-up solution} and Section \ref{Extinction of the minimal blow-up solution}. Particularly, we give the precise construction of the MEI-functional $\mD$, which will help us to set up the inductive hypothesis given in Section \ref{Existence of the minimal blow-up solution}.
\begin{lemma}\label{positive or negative k}
Let $u\in H^1(\R^2)\setminus\{0\}$ with $\mM(u)<{\mM(Q)}$. Then there exists a unique $\ld(u)>0$ such that
\begin{equation}
\mK(T_\ld u)\left\{
             \begin{array}{ll}
             >0, &\text{if $\ld\in(0,\ld(u))$},  \\
             =0,&\text{if $\ld=\ld(u)$},\\
             <0,&\text{if $\ld\in(\ld(u),\infty)$}.
             \end{array}
\right.
\end{equation}
\end{lemma}

\begin{proof}
We first obtain that
\begin{align}
\mK(T_\ld u)&=\ld^2 \|\nabla u\|_2^2-\frac{\ld^2}{2} \|u\|_\tas^\tas-\frac{2\ld^4}{3}\|u\|_\tbs^\tbs.\label{one a}
\end{align}
By \eqref{GN-L1} we have
\begin{align}\label{GN-L2}
2\|\nabla u\|_2^2-\|u\|_{\tas}^\tas\geq 2\bg(1-\frac{\mM(u)}{{\mM(Q)}}\bg)\|\nabla u\|_2^2>0.
\end{align}
Thus one easily sees that $\mK(T_\ld u)$ is positive on $(0,\ld(u))$ and negative on $(\ld(u),\infty)$, where
\begin{align}
\ld(u)=\bg(\frac{3(2\|\nabla u\|_2^2-\|u\|_{\tas}^\tas)}{4\|u\|_6^6}\bg)^{\frac{1}{2}}
\end{align}
is the unique zero of $\ld\mapsto \mK(T_\ld u)$ on $(0,\infty)$. This completes the proof.
\end{proof}

\begin{lemma}\label{pos of k implies pos of h}
Assume that $\mK(u)\geq 0$. Then $\mH(u)\geq 0$. If additionally $\mK(u)> 0$, then also $\mH(u)> 0$.
\end{lemma}
\begin{proof}
We have
\begin{align}
\mH(u)\geq \mH(u)-\frac{1}{2}\mK(u)=\frac{1}{6}\|u\|_{\tbs}^\tbs\geq 0.
\end{align}
It is straightforward to obtain that the last inequality can be replaced by the strict one when $u\neq 0$, which is the case when $\mK(u)>0$.
\end{proof}

\begin{lemma}\label{bound of gradient by energy}
Let $\delta\in(0,1)$ and let $u\in \mA$. Suppose also that
\begin{align}
\mM(u)\leq (1-\delta)\mM(Q)
\end{align}
with some $\delta\in(0,1)$. Then
\begin{align}
\|u\|_\tbs^\tbs&< \frac{3}{2}\|\nabla u\|_2^2,\label{bound of tbs by 2}\\
\|u\|_\tas^\tas&\leq 2(1-\delta)\|\nabla u\|_2^2,\\
\frac{\delta}{4} \|\nabla u\|_2^2&< \mH(u)\leq \frac{1}{2}\|\nabla u\|_2^2.\label{third}
\end{align}
\end{lemma}

\begin{proof}
The first inequality follows immediately from the fact that $\mK(u)>0$ for $u\in{\mA}$. For the second one, we obtain that
\begin{align}
&\|u\|_\tas^\tas\nonumber\\
\leq&\, \mathrm{C}^{-1}_{\mathrm{GN}}\|\nabla u\|_2^2\|u\|_2^{2}\nonumber\\
\leq& \,2{\mM(Q)}^{-1}\|\nabla u\|_2^2(1-\delta){\mM(Q)}\nonumber\\
=&\,2(1-\delta)\|\nabla u\|_2^2.
\end{align}
The first inequality in \eqref{third} follows from
\begin{align}
\mH(u)&> \mH(u)-\frac{1}{4}\mK(u)\nonumber\\
&=\frac{1}{8}(2\|\nabla u\|_2^2-\|u\|_\tas^\tas)\nonumber\\
&\geq \frac{1}{4}\bg(1-\frac{\mM(u)}{{\mM(Q)}}\bg)\|\nabla u\|_2^2\nonumber\\
&\geq \frac{\delta}{4}\|\nabla u\|_2^2
\end{align}
and the second inequality in \eqref{third} follows immediately from the non-positivity of the cubic and quintic nonlinearities.
\end{proof}

\begin{lemma}\label{monotone lemma}
The mapping $c\mapsto m_c$ is continuous and monotone decreasing on $(0,{\mM(Q)})$.
\end{lemma}

\begin{proof}
The proof follows the arguments of \cite{Bellazzini2013}, where we also need to take the effect of the mass constraint into account. We first show that the function $f$ defined by
\begin{align*}
f(a,b):=\max_{t>0}\{at^2-bt^{4}\}
\end{align*}
is continuous on $(0,\infty)^2$. Define
\begin{align*}
g(a,b,t):=at^2-bt^{4}.
\end{align*}
Then for any $a,b>0$, there exists a unique $t_*>0$ such that
\begin{align}
\pt_t g(t_*,a,b)&=0,\\
\pt_{tt}g(t_*,a,b)&<0.
\end{align}
By the implicit function theorem we deduce the existence of a continuous function $h$ in a neighborhood of $(a,b)$ such that $\pt_t g(h(a,b),a,b)=0$. Hence
\begin{align}
f(a,b)=g(h(a,b),a,b)
\end{align}
and therefore the mapping $(a,b)\mapsto f(a,b)$ is continuous. Next we show that for any $0<c_1<c_2<{\mM(Q)}$ and $\vare>0$ we have
\begin{align}
m_{c_2}\leq m_{c_1}+\vare.
\end{align}
Define the set $V(c)$ by
\begin{align*}
V(c):=\{u\in H^1(\R^2):\mM(u)=c,\mK(u)=0\}.
\end{align*}
By the definition of $m_{c_1}$ there exists some $u_1\in V(c_1)$ such that
\begin{align}\label{pert 2}
\mH(u_1)\leq m_{c_1}+\frac{\vare}{2}.
\end{align}
Let $\eta\in C^{\infty}_0(\R^2;[0,1])$ be a cut-off function with $\eta(x)=1$ for $|x|\leq 1$ and $\eta(x)=0$ for $|x|\geq 2$. For $\delta>0$, define
\begin{equation*}
\tilde{u}_{1,\delta}(x):= \eta(\delta x)\cdot u_1(x).
\end{equation*}
Then $\tilde{u}_{1,\delta}\to u_1$ in $H^1(\R^2)$ as $\delta\to 0$. Therefore,
\begin{align}
\|\nabla\tilde{u}_{1,\delta}\|_2^2&\to \|\nabla u_1\|_2^2, \\
\|\tilde{u}_{1,\delta}\|_p&\to \| u_1\|_p
\end{align}
for all $p\in[2,\tbs]$ as $\delta\to 0$. Using \eqref{GN-L2} we know that
\begin{align}
\frac{1}{2}\|\nabla v\|_2^2>\frac{1}{\tas}\|v\|_\tas^\tas
\end{align}
for all $v\in H^1(\R^2)$ with $\mM(v)\in(0,{\mM(Q)})$. Since $c_1\in(0,{\mM(Q)})$, we infer that $\mM(\tilde{u}_{1,\delta})\in (0,{\mM(Q)})$ for sufficiently small $\delta$. Combining with the continuity of the function $f$ given previously we conclude that
\begin{align}\label{pert 1}
\max_{t>0}\mH(T_t\tilde{u}_{1,\delta})&=\max_{t>0}\{t^2(\frac{1}{2}\|\nabla\tilde{u}_{1,\delta}\|_2^2-\frac{1}{\tas}\|\tilde{u}_{1,\delta}\|_\tas^\tas)
-\frac{t^4}{\tbs}\|\tilde{u}_{1,\delta}\|_\tbs^\tbs\}\nonumber\\
&\leq \max_{t>0}\{t^2(\frac{1}{2}\|\nabla u_1\|_2^2-\frac{1}{\tas}\|u_1\|_\tas^\tas)
-\frac{t^4}{\tbs}\|u_1\|_\tbs^\tbs\}+\frac{\vare}{4}\nonumber\\
&=\max_{t>0}\mH(T_t u_1)+\frac{\vare}{4}
\end{align}
for sufficiently small $\delta>0$. Now let $v\in C_0^\infty(\R^2)$ with $\mathrm{supp}(v)\subset \R^2\backslash B(0,2^{-1}\delta)$ and define
\begin{equation*}
v_0:= \bg(\frac{c_2-\mM(\tilde{u}_{1,\delta})}{\mM(v)}\bg)^{\frac12}\,v.
\end{equation*}
We have $\mM(v_0)=c_2-\mM(\tilde{u}_{1,\delta})$. Let
\begin{align*}
w_\ld:=\tilde{u}_{1,\delta}+T_\ld v_0
\end{align*}
with some to be determined $\ld\in(0,1)$. By definition one easily sees that for all $\ld\in(0,1)$ the supports of $\tilde{u}_{1,\delta}$ and $T_\ld v_0$ are disjoint, thus
\begin{align}
\|w_\ld\|_p=\|\tilde{u}_{1,\delta}\|_p+\|T_\ld v_0\|_p
\end{align}
for all $p\in [2,6]$. Particularly we infer that $\mM(w_\ld)=c_2$. Moreover one easily verifies that
\begin{align}
\|\nabla w_\ld\|_2&\to\|\nabla \tilde{u}_{1,\delta}\|_2,\\
\| w_\ld\|_p&\to\| \tilde{u}_{1,\delta}\|_p
\end{align}
for all $p\in(2,6]$ as $\ld\to 0$. Using the continuity of the function $f$ once again we obtain that
\begin{align}
\max_{t>0}\mH(T_tw_\ld)\leq \max_{t>0}\mH(T_t \tilde{u}_{1,\delta})+\frac{\vare}{4}
\end{align}
for sufficiently small $\ld>0$. Finally, combing with \eqref{pert 2} and \eqref{pert 1} we conclude that
\begin{align}
m_{c_2}&\leq \max_{t>0}\mH(T_t w_\lambda)\leq \max_{t>0}\mH(T_t\tilde{u}_{1,\delta})+\frac{\varepsilon}{4}\nonumber\\
&\leq \max_{t>0}\mH(T_tu_1)+\frac{\varepsilon}{2}=\mH(u_1)+\frac{\varepsilon}{2}\leq m_{c_1}+\varepsilon.
\end{align}
Choosing $\vare$ arbitrarily small then completes the proof of the monotonicity. The arguments for proving the continuity of the mapping $c\mapsto m_c$ are very similar to the previous ones, we therefore omit the details of the straightforward but tedious modification and refer for instance to \cite[Lem. 5.4]{Bellazzini2013} or \cite[Lem. 3.3]{SoaveSubcritical} for a complete proof.
\end{proof}

The following lemma shows that the NLS-flow leaves solutions starting from ${\mA}$ invariant.
\begin{lemma}\label{invariant lemma}
Let $u$ be a solution of \eqref{NLS} such that $u(0)\in {\mA}$. Then $u(t)\in {\mA}$ for all $t$ in the maximal lifespan $I_{\max}$. Assume also $\mM(u)=(1-\delta)\mM(Q)$. Then
\begin{align}\label{lower bound kt}
\mK(u(t))\geq\min\bg\{\delta\mH(u(0)),\,\bg(\bg(\frac{2}{\delta}\bg)^{\fraco}-1\bg)^{-1}\bg(m_{\mM(u(0))}-\mH(u(0))\bg)\bg\}.
\end{align}
for all $t\in I_{\max}$.
\end{lemma}

\begin{proof}
By the mass and energy conservation, to show the invariance of $\mA$ under the NLS-flow we only need to show that $\mK(u(t))> 0$ for all $t\in I_{\max}$. Suppose that there exists some $t\in I_{\max}$ such that $\mK(u(t))\leq 0$. By continuity of $u(t)$ there exists some $s\in(0,t]$ such that $\mK(u(s))=0$. By conservation of mass we also know that $0<\mM(u(s))<{\mM(Q)}$. Now using the definition of $m_c$ we immediately obtain that
\begin{align}
m_{\mM(u(s))}\leq \mH(u(s))<m_{\mM(u(0))}= m_{\mM(u(s))},
\end{align}
a contradiction. We now show \eqref{lower bound kt}. Direct calculation yields
\begin{align}\label{der of h}
\frac{d^2}{d\ld^2}\mH(T_\ld u(t))=-\frac{1}{\ld^2}\mK(T_\ld u(t))+\frac{2}{\ld^2}\bg(\mK(T_\ld u(t))-\frac{2}{3}\|T_\ld u(t)\|_\tbs^\tbs\bg).
\end{align}
If
\begin{align}
\mK(u(t))-\frac{2}{3}\| u(t)\|_\tbs^\tbs\geq 0,
\end{align}
then using \eqref{GN-L2} we deduce that
\begin{align}
\mK(u(t))&=\|\nabla u\|_2^2-\frac{1}{2}\|u\|_\tas^\tas-\frac{2}{3}\|u\|_\tbs^\tbs\nonumber\\
&\geq \delta\|\nabla u\|_2^2-\mK(u(t)),
\end{align}
which combining with \eqref{third} implies
\begin{align}\label{lower bd 1}
\mK(u(t))\geq \frac{\delta}{2}\|\nabla u(t)\|_2^2\geq  \delta\mH(u(0)),
\end{align}
where for the last inequality we also used the conservation of energy. Suppose now that
\begin{align}\label{negative 0}
\mK(u(t))-\frac{2}{3}\| u(t)\|_\tbs^\tbs< 0.
\end{align}
Then
\begin{align}
\frac{2}{3}\| u(t)\|_\tbs^\tbs&>\|\nabla u(t)\|_2^2-\frac{1}{2}\|u\|_\tas^\tas-\frac{2}{3}\|u\|_\tbs^\tbs\nonumber\\
&\geq \delta\|\nabla u(t)\|_2^2-\frac{2}{3}\|u(t)\|_\tbs^\tbs,
\end{align}
hence
\begin{align}\label{2 by 2*}
\|u(t)\|_\tbs^\tbs> \frac{3\delta}{4}\|\nabla u(t)\|_2^2.
\end{align}
Since $\mK(u(t))> 0$, by Lemma \ref{positive or negative k} we know that there exists some $\ld_*\in(1,\infty)$ such that
\begin{align}\label{positive ld for K}
\mK(T_\ld u(t))\geq 0\quad\forall\,\ld\in[1,\ld_*]
\end{align}
and
\begin{align}
0&=\mK(T_{\ld_*}u(t))\nonumber\\
&=\ld_*^2(\|\nabla u(t)\|_2^2-\frac{1}{2}\|u(t)\|_\tas^\tas)-\frac{2\ld_*^4}{3}\|u(t)\|_\tbs^\tbs,
\end{align}
which in turn gives
\begin{align}\label{2* by 2}
\|u(t)\|_\tbs^\tbs&= \frac{3\ld_*^{-2}}{2}(\|\nabla u(t)\|_2^2-\frac{1}{2}\|u(t)\|_\tas^\tas)\leq \frac{3\ld_*^{-2}}{2}\|\nabla u(t)\|_2^2.
\end{align}
\eqref{2 by 2*} and \eqref{2* by 2} then result in
\begin{align}\label{bound of ld}
\ld_*< \bg(\frac{2}{\delta}\bg)^{\fraco}.
\end{align}
On the other hand, direct calculation yields
\begin{align}\label{negative 1}
\frac{d}{ds}\bg(\frac{1}{s^2}\bg(\mK(T_s u(t))-\frac{2}{3}\|T_s u(t)\|_\tbs^\tbs\bg)\bg)=-\frac{8}{3}s \|u(t)\|_\tbs^\tbs<0
\end{align}
for $s>0$. Integrating \eqref{negative 1} over $[1,\ld]$ and using \eqref{negative 0}, we find that for $\ld\geq 1$
\begin{align}\label{negative der of k}
\frac{1}{\ld^2}\bg(\mK(T_\ld u(t))-\frac{2}{3}\|T_\ld u(t)\|_\tbs^\tbs\bg)\leq 0.
\end{align}
\eqref{der of h}, \eqref{positive ld for K} and \eqref{negative der of k} imply that $\frac{d^2}{d\ld^2}\mH(T_\ld u(t))\leq 0$ for all $\ld\in[1,\ld_*]$. Finally, combining with \eqref{bound of ld}, the fact that $\mK(T_{\ld_*}u(t))=0$ and Taylor expansion we infer that
\begin{align}
&\,\bg(\bg(\frac{2}{\delta}\bg)^{\fraco}-1\bg)\mK(u(t))\nonumber\\
\geq&\, \,\bg(\ld_*-1\bg)\bg(\frac{d}{d\ld}|_{\ld=1}\mH(T_\ld u(t))\bg)\nonumber\\
\geq&\,\, \mH(T_{\ld_*}u(t))-\mH(u(t))\nonumber\\
\geq&\,\, m_{\mM(u(0))}-\mH(u(0)).
\end{align}
This together with \eqref{lower bd 1} yields \eqref{lower bound kt}.
\end{proof}

\begin{lemma}\label{mtilde equal m}
Let
\begin{align}
\tilde{m}_{c}&:=\inf_{u\in H^1(\R^2)}\{\mI(u):\|u\|_2= c,\mK(u)\leq 0\}.
\end{align}
Then $m_{c}=\tm_c$.
\end{lemma}

\begin{proof}
Let $(u_n)_n$ be a minimizing sequence for the variational problem of $\tilde{m}_{\delta}$, i.e.
\begin{align}
\lim_{n\to\infty}\mI(u_n)&=\tilde{m}_c,\\
\mM(u_n)&=c\quad\forall\,n\in\N,\\
\mK(u_n)&\leq 0\quad\forall\,n\in\N.
\end{align}
Using Lemma \ref{positive or negative k} we know that there exists some $\ld_n\in(0,1]$ such that $\mK(T_{\ld_n}u_n)$ is equal to zero. Thus
\begin{align}
m_c\leq \mH(T_{\ld_n}u_n)=\mI(T_{\ld_n}u_n)\leq \mI(u_n)=\tilde{m}_c+o_n(1).
\end{align}
Sending $n\to\infty$ we infer that $m_c\leq \tm_c$. On the other hand,
\begin{align}
\tm_c
&\leq\inf_{u\in H^1(\R^2)}\{\mI(u):\mM(u)=c,\mK(u)=0\}\nonumber\\
&=\inf_{u\in H^1(\R^2)}\{\mH(u):\mM(u)=c,\mK(u)=0\}=m_c.
\end{align}
This completes the proof.
\end{proof}

We now define
\begin{align*}
m_0&:=\lim_{c\to 0}m_c\in(0,\infty],\\
m_{Q}&:=\lim_{c\to{\mM(Q)}}m_c\in[0,\infty).
\end{align*}
Also define the set $\Omega$ by its complement
\begin{align*}
\Omega^c&:=\{(c,h)\in\R^2: c\geq {\mM(Q)}\}\cup\{(c,h)\in\R^2: c\geq 0 \wedge h\geq m_Q\}.
\end{align*}
The MEI-functional $\mD:\R^2\to\R\cup\{\infty\}$ is defined by
\begin{align*}
\mD(c,h)=\left\{
             \begin{array}{ll}
             h+\frac{h+c}{\mathrm{dist}((c,h),\Omega^c)},&\text{if $(c,h)\in \Omega$},\\
             \infty,&\text{otherwise}.
             \end{array}
\right.
\end{align*}
For $u\in H^1(\R^2)$ we simply write $\mD(u):=\mD(\mM(u),\mH(u))$. A schematic description of the domain $\Omega$ is given by Fig. \ref{dihe} below.

\begin{figure}[ht!]
  \centering
  \includegraphics[width=30mm]{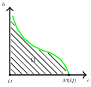}
  \caption{The green curve represents the mapping $c\mapsto m_c$ and the shadow region is the intersection of $\Omega$ and $(0,\infty)^2$, which corresponds to the set $\mA$ due to Lemma \ref{killip visan curve} (ii).}\label{dihe}
\end{figure}

\begin{remark}
By modifying the arguments in \cite[Lem. 5.6]{Bellazzini2013} and \cite[Lem. 3.3]{wei2021normalized} one is able to show
\begin{align}
m_0=\infty,\quad m_{Q}=0.
\end{align}
Nevertheless, the precise values of $m_0$ and $m_Q$ have no impact on the scattering result; All we need here is the monotonicity and continuity of the curve $c\mapsto m_c$. We will therefore postpone the proof to \ref{m0 mQ}.
\end{remark}


We end this section by establishing some useful properties of the MEI-functional $\mD$.
\begin{lemma}\label{killip visan curve}
Suppose that $v\in H^1(\R^2)$ and satisfies $\mK(v)\geq 0$. Then
\begin{itemize}
\item[(i)]$\mD(v)=0$ if and only if $u=0$.
\item[(ii)] $0<\mD(v)<\infty$ if and only if $v\in\mA$.
\item[(iii)] If $u$ is a solution of \eqref{NLS} with $u(t_0)\in\mA$, then $\mD(u(t))=\mD(u(t_0))$ for all $t\in I_{\max}$.
\item[(iv)] Let $u_1,u_2\in\mA$ with $\mM(u_1)\leq \mM(u_2)$ and $\mH(u_1)\leq \mH(u_2)$. Then $\mD(u_1)\leq \mD(u_2)$. If additionally either $\mM(u_1)<\mM(u_2)$ or $\mH(u_1)<\mH(u_2)$, then $\mD(u_1)<\mD(u_2)$.
\item[(v)] Let $\mD_0\in(0,\infty)$. Then
\begin{align}
\|\nabla u\|^2_{2}&\sim_{\mD_0}\mH(u),\\
\|u\|^2_{H^1}&\sim_{\mD_0}\mH(u)+\mM(u)\sim_{\mD_0}\mD(u)
\end{align}
uniformly for all $u\in \mA$ with $\mD(u)\leq \mD_0$.
\item[(vi)] Let $\mD_0\in(0,\infty)$. Then
\begin{align}\label{small of unaaa}
m_{\mM(u)}-\mH(u)\gtrsim_{\mD_0} 1
\end{align}
uniformly for all $u\in \mA$ with $\mD(u)\leq \mD_0$.
\end{itemize}
\end{lemma}

\begin{proof}
\begin{itemize}
\item[(i)]It is trivial that $v=0$ implies $\mD(v)=0$. Suppose now $\mD(v)=0$. Since $\mK(v)\geq 0$, we infer from Lemma \ref{pos of k implies pos of h} that $\mH(v)\geq 0$. In this case, $\mD(v)=0$ can only happen when $v=0$.
\item[(ii)]It is trivial that $v\in \mA$ implies $\mD(v)<\infty$. By Lemma \ref{pos of k implies pos of h} we also know that $\mH(v)> 0$ for $v\in\mA$, which implies $\mD(v)> 0$ by the definition of $\mD$. Now let $\0<\mD(v)<\infty$. Then $\mM(v)\in(0,{\mM(Q)})$. By definition of $\mD$ we also know that $\mH(v)<m_{\mM(v)}$. Using the definition of $m_{\mM(v)}$ we conclude that $\mK(v)>0$ (under the precondition $\mK(v)\geq 0$) and therefore $v\in\mA$.
\item[(iii)]This follows immediately from the conservation of mass and energy of the NLS flow, the definition of $\mD$ and Lemma \ref{invariant lemma}.
\item[(iv)]This follows from the fact that $c\mapsto m_c$ is monotone decreasing on $(0,{\mM(Q)})$ and the definition of $\mD$.
\item[(v)]Since $u\in\mA$, we know that $\mM(u)\in(0,\mM(Q))$ and using Lemma \ref{pos of k implies pos of h} also $\mH(u)\in (0,m_{\mM(u)})$. Thus
\begin{align}
&\,\mathrm{dist}\bg((\mM(u),\mH(u)),\Omega^c\bg)\nonumber\\
\leq&\, \mathrm{dist}\bg((\mM(u),\mH(u)),(\mM(Q),\mH(u))\bg)=\mM(Q)-\mM(u).
\end{align}
Since $\mH(u)\geq 0$, we have
\begin{align}\label{mass constraint}
\mD(u)\geq \frac{\mM(u)}{\mM(Q)-\mM(u)},
\end{align}
which implies
\begin{align}
\frac{1}{\mD(u)+1}\leq \frac{\mM(Q)-\mM(u)}{\mM(Q)}=1-\frac{\mM(u)}{\mM(Q)}.
\end{align}
Since $\mK(u)> 0$, we have
\begin{align}\label{energy constraint}
\mD(u)\geq& \,\mH(u)> \mH(u)-\frac{1}{4}\mK(u)\nonumber\\
=&\,\frac{1}{8}(2\|\nabla u\|_2^2-\|u\|_\tas^\tas)\nonumber\\
\geq &\,\frac{1}{4}\bg(1-\frac{\mM(u)}{{\mM(Q)}}\bg)\|\nabla u\|_2^2\nonumber\\
\geq &\,\frac{\|\nabla u\|_2^2}{4(\mD(u)+1)},
\end{align}
therefore $\|\nabla u\|_2^2\lesssim_{\mD_0}\mH(u)$. Combining with \eqref{third} and (iii) we have
\begin{align}
\|\nabla u\|_2^2&\sim_{\mD_0}\mH(u),\\
\|u\|_{H^1}^2&\sim_{\mD_0}\mH(u)+\mM(u).
\end{align}
It remains to show $\mH(u)+\mM(u)\sim_{\mD_0}\mD(u)$. Using \eqref{mass constraint} and \eqref{energy constraint} we infer that
\begin{align}
\mH(u)+\mM(u)\sim_{\mD_0}\|u\|^2_{H^1}\lesssim_{\mD_0} \mD(u).
\end{align}
To show $\mD(u)\lesssim_{\mD_0} \mH(u)+\mM(u)$ we discuss the following different cases: If $\mM(u)\geq\frac{1}{2}\mM(Q)$, then using the fact that $\mH(u)\geq 0$ we have
\begin{align}
\mathrm{dist}\bg((\mM(u),\mH(u)),\Omega^c\bg)\geq \frac{\mM(u)}{\mD_0}\geq\frac{\mM(Q)}{2\mD_0},
\end{align}
which implies
\begin{align}
\mD(u)\leq \frac{2\mD_0}{\mM(Q)}\bg(\mM(u)+\mH(u)\bg)+\mH(u).
\end{align}
If $\mM(u)\leq \frac{1}{2}\mM(Q)$ and $\mH(u)\geq \frac{1}{2}m_{\frac{1}{2}\mM(Q)}$, then analogously we obtain
\begin{align}
\mD(u)\leq \frac{2\mD_0 }{m_{\frac{1}{2}\mM(Q)}}\bg(\mM(u)+\mH(u)\bg)+\mH(u).
\end{align}
If $\mM(u)\leq \frac{1}{2}\mM(Q)$ and $\mH(u)\leq \frac{1}{2}m_{\frac{1}{2}\mM(Q)}$, then using the monotonicity of $c\mapsto m_c$ we have
\begin{align}
\mathrm{dist}\bg((\mM(u),\mH(u)),\Omega^c\bg)\geq\mathrm{dist}\bg(\bg(\frac{1}{2}\mM(Q),\frac{1}{2}m_{\frac{1}{2}\mM(Q)}\bg),\Omega^c\bg)=:\alpha_0>0.
\end{align}
Therefore
\begin{align}
\mD(u)\leq \frac{1}{\alpha_0}\bg(\mM(u)+\mH(u)\bg)+\mH(u),
\end{align}
which completes the proof of (v).
\item[(vi)]
If this were not the case, then we could find a sequence $(u_n)_n\subset\mA$ such that
\begin{align}\label{small of un}
m_{\mM(u_n)}-\mH(u_n)=o_n(1),
\end{align}
which implies
\begin{align}
\mathrm{dist}\bg(\bg(\mM(u_n),\mH(u_n)\bg),\Omega^c\bg)\leq&\, \mathrm{dist}\bg(\bg(\mM(u_n),\mH(u_n)\bg),\bg(\mM(u_n),m_{\mM(u_n)}\bg)\bg)\nonumber\\
=&\,m_{\mM(u_n)}-\mH(u_n)=o_n(1).
\end{align}
If $\mM(u_n)\gtrsim 1$, then
\begin{align}\label{on1 contra}
\mD(u_n)\gtrsim \frac{1}{o_n(1)},
\end{align}
contradicting $\mD(u_n)\leq \mD_0$. If $\mM(u_n)=o_n(1)$, then by the monotonicity of $c\mapsto m_c$ and \eqref{small of un} we know that $\mH(u_n)\gtrsim 1$ and similarly we may again derive the contradiction \eqref{on1 contra}.
\end{itemize}
\end{proof}

\section{Existence of the minimal blow-up solution}\label{Existence of the minimal blow-up solution}
We define
\begin{align}\label{def of tau}
\tau(\mD_0):=\sup\bg\{&\|\la\nabla\ra^{\frac{1}{2}} \psi\|_{{L_{t,x}^4}(I_{\max})}:\nonumber\\
&\quad\quad\quad\quad\text{$\psi$ is solution of \eqref{NLS}, }\psi(0)\in {\mA},\mD(\psi(0))\leq \mD_0\bg\}
\end{align}
and
\begin{align}\label{introductive hypothesis}
\mD^*&:=\sup\{\mD_0>0:\tau(\mD_0)<\infty\}.
\end{align}
By Lemma \ref{well posedness lemma}, Remark \ref{remark} and Lemma \ref{killip visan curve} (v) we know that $\tau(\mD_0)<\infty$ for sufficiently small $\mD_0$. We therefore assume that $\mD^*<\infty$ and aim to derive a contradiction, which will imply $\mD^*=\infty$ and the desired proof is complete in view of Lemma \ref{killip visan curve} (ii). By the inductive hypothesis we may find a sequence $(\psi_n)_n$ with $(\psi_n(0))_n\subset {\mA}$ which are solutions of \eqref{NLS} with maximal lifespan $(I_{n})_n$ such that
\begin{gather}
\lim_{n\to\infty}\|\la\nabla\ra^{\frac{1}{2}}\psi_n\|_{{L_{t,x}^4}(I_n)}=\infty,\label{oo1}\\
\lim_{n\to\infty}\mD(\psi_n(0))=\mD^*.\label{oo2}
\end{gather}
Up to a subsequence we may also assume that
\begin{align}\label{convergence}
(\mM(\psi_n(0)),\mH(\psi_n(0)))\to(\mM_0,\mH_0)\quad\text{as $n\to\infty$}.
\end{align}
By continuity of $\mD$ and finiteness of $\mD^*$ we know that
\begin{align}
\mD^*&=\mD(\mM_0,\mH_0),\\
\mM_0&\in(0,\mM(Q)),\\
\mH_0&\in[0,m_{\mM_0}).
\end{align}
By Lemma \ref{killip visan curve} (v) we deduce that $(\psi_n(0))_n$ is a bounded sequence in $H^1(\R^2)$. Using Lemma \ref{linear profile} applied to $(\psi_n(0))_n$ we infer that there exist some number $K^*\in\N\cup\{\infty\}$, a sequence of nonzero linear profiles $(\phi^j)_{1\leq j\leq K^*}\subset L^2(\R^2)$, a sequence of symmetry parameters $(\ld_n^j,t_n^j,x_n^j,\xi_n^j)_{n\in\N,1\leq j\leq K^*}\subset(0,\infty)\times\R\times \R^2\times\R^2$ with $\sup_{n\in\N}|\xi_n^j|\lesssim_j 1$ and a sequence of remainders $(w_{n}^k)_{n\in\N,1\leq k\leq K^*}\subset H^1(\R^2)$ such that
\begin{itemize}
\item[(i)]The parameters $(\ld_n^j,t_n^j,x_n^j,\xi_n^j)_{n,j}$ satisfy
\begin{align}\label{orthog of pairs}
\lim_{n\to\infty}&\bg\{\bg|\log \frac{\ld_n^j}{\ld_n^l}\bg|+\frac{|t_n^j-t_n^l|}{(\ld_n^j)^2}\nonumber\\
&\quad\quad\quad+\ld_n^j|\xi_n^j-\xi_n^l|+\frac{|x_n^j-x_n^l+2t_n^j(\xi_n^j-\xi_n^l)|}{\ld_n^j}\bg\}=\infty
\end{align}
for all finite $1\leq j,l\leq K^*$ with $j\neq l$. Moreover,
\begin{align}
\lim_{n\to\infty}\ld_n^j&=\ld^j_\infty\in\{1,\infty\},\\
\ld_n^j&\equiv1\quad\text{if $\ld_\infty^j=1$}.
\end{align}

\item[(ii)] There exists some $\theta\in(0,1)$ such that for any finite $1\leq k\leq K^*$ we have the decomposition
\begin{align}
\psi_n(0)=\sum_{j=1}^k T_n^j P_n^j \phi^j+w_n^k,
\end{align}
where
\begin{align}
T_n^j u(x):=e^{ix\cdot\xi_n^j}e^{-t_n^j\Delta}(\ld_n^j)^{-1}u\bg((\ld_n^j)^{-1}(\cdot-x_n^j)\bg)(x)
\end{align}
and
\begin{equation}
P_n^j\phi^j=\left\{
             \begin{array}{ll}
             \phi^j,&\text{if $\ld_\infty^j=1$},\\
             P_{\leq(\ld_n^j)^\theta}\phi^j,&\text{if $\ld_\infty^j=\infty$}.
             \end{array}
\right.
\end{equation}
Moreover, if $\ld^j_\infty=1$, then $\xi_n^j\equiv 0$ and $\phi^j\in H^1(\R^2)$.
\item[(iii)] The remainders $(w_n^k)_{n,k}$ satisfy
\begin{align}\label{to zero wnk}
\lim_{k\to K^*}\lim_{n\to\infty}\|\la\nabla\ra^{\frac{1}{2}} e^{it\Delta}w_n^k\|_{{L_{t,x}^4}(\R)}=0.
\end{align}

\item[(iv)] The following orthogonal properties are satisfied: for $s\in\{0,1\}$ and finite $1\leq J\leq K^*$ we have
\begin{align}
\||\nabla|^s\psi_n\|_2^2&=\sum_{j=1}^J \||\nabla|^s T_n^j P_n^j\phi^j\|_2^2+\||\nabla |^s w_n^J\|_2^2+o_n(1),\label{orthog L2 and H1}\\
\mH(\psi_n)&=\sum_{j=1}^J \mH(T_n^j P_n^j\phi^j)+\mH( w_n^J)+o_n(1),\label{conv of h}\\
\mI(\psi_n)&=\sum_{j=1}^J \mI(T_n^j P_n^j\phi^j)+\mI( w_n^J)+o_n(1),\label{conv of i}\\
\mK(\psi_n)&=\sum_{j=1}^J \mK(T_n^j P_n^j\phi^j)+\mK( w_n^J)+o_n(1).\label{conv of k}
\end{align}

\end{itemize}
\begin{lemma}\label{category 0 and 1}
There exists a global solution $u_c$ of \eqref{NLS} such that
\begin{align}
(\mD(u_c),\mM(u_c),\mH(u_c))=(\mD^*,\mM_0,\mH_0)\label{requirement1}
\end{align}
and
\begin{align}
\|\la\nabla\ra^{\frac{1}{2}} u_c\|_{{L_{t,x}^4}(\R)}=\infty.\label{requirement2}
\end{align}
\end{lemma}

\begin{proof}
We first show that for a given nonzero linear profile $\phi^j$ we have
\begin{align}
\mH(T_n^jP_n^j\phi^j)&> 0,\label{bd for S}\\
\mK(T_n^jP_n^j\phi^j)&> 0\label{pos of K}
\end{align}
for all sufficiently large $n=n(j)\in\N$. Since $\phi^j\neq 0$ we know that $T_n^jP_n^j\phi^j\neq 0$ for sufficiently large $n$. Suppose now that \eqref{pos of K} does not hold. Up to a subsequence we may assume that $\mK(T_n^jP_n^j\phi^j)\leq 0$ for all sufficiently large $n$. By the non-negativity of $\mI$, \eqref{conv of i} and \eqref{small of unaaa} we know that there exists some sufficiently small $\delta>0$ depending on $\mD^*$ and some sufficiently large $N_1$ such that for all $n>N_1$ we have
\begin{align}\label{contradiction1}
\tm_{\mM(T_n^jP_n^j\phi^j)}&\leq\mI(T_n^jP_n^j\phi^j)\leq \mI(\psi_n(0))+\delta\nonumber\\
&\leq\mH(\psi_n(0))+\delta\leq m_{\mM(\psi_n(0))}-2\delta,
\end{align}
where $\tm$ is the quantity defined by Lemma \ref{mtilde equal m}. By continuity of $c\mapsto m_c$ we also know that for sufficiently large $n$ we have
\begin{align}\label{contradiction3}
m_{\mM(\psi_n(0))}-2\delta\leq m_{\mM_0}-\delta.
\end{align}
Using \eqref{orthog L2 and H1} we deduce that for any $\vare>0$ there exists some large $N_2$ such that for all $n>N_2$ we have
\begin{align}
\mM(T_n^jP_n^j\phi^j)\leq \mM_0+\vare.
\end{align}
From the continuity and monotonicity of $c\mapsto m_c$ and Lemma \ref{mtilde equal m}, we may choose some sufficiently small $\vare$ to see that
\begin{align}\label{contradiction2}
\tm_{\mM(T_n^jP_n^j\phi^j)}=m_{\mM(T_n^jP_n^j\phi^j)}\geq m_{\mM_0+\vare}\geq m_{\mM_0}-\frac{\delta}{2}.
\end{align}
Now \eqref{contradiction1}, \eqref{contradiction3} and \eqref{contradiction2} yield a contradiction. Thus \eqref{pos of K} holds, which combining with Lemma \ref{pos of k implies pos of h} also yields \eqref{bd for S}. Similarly, for each $j\in\N$ we have
\begin{align}
\mH(w_n^j)&\geq 0,\label{bd for S wnj} \\
\mK(w_n^j)&\geq 0\label{pos of K wnj}
\end{align}
for sufficiently large $n$. Now using \eqref{convergence} we have for any $k\in\N$
\begin{align}
\mM_0&=\sum_{j=1}^k \mM(T_n^jP_n^j\phi^j)+\mM(w_n^k)+o_n(1),\label{mo sum}\\
\mH_0&=\sum_{j=1}^k \mH(T_n^jP_n^j\phi^j)+\mH(w_n^k)+o_n(1).\label{eo sum}
\end{align}
From \eqref{mo sum} and \eqref{eo sum} we infer that two different scenarios will potentially take place: either
\begin{align}
\sup_{j\in\N}\lim_{n\to\infty}\mM(T_n^jP_n^j\phi^j)&=\mM_0\text{ and}\nonumber\\
\sup_{j\in\N}\lim_{n\to\infty}\mH(T_n^jP_n^j\phi^j)&=\mH_0,\label{first situation}
\end{align}
or there exists some $\delta>0$ such that
\begin{align}
\sup_{j\in\N}\lim_{n\to\infty}\mM(T_n^jP_n^j\phi^j)&\leq \mM_0-\delta\text{ or}\nonumber\\
\sup_{j\in\N}\lim_{n\to\infty}\mH(T_n^jP_n^j\phi^j)&\leq\mH_0-\delta.\label{second situation}
\end{align}
We show that starting from \eqref{first situation} one is able to derive a minimal blow-up solution $u_c$ which satisfies \eqref{requirement1} and \eqref{requirement2}, while from \eqref{second situation} we get a contradiction. We begin firstly with \eqref{first situation}. In this case, since the summands in \eqref{mo sum} and \eqref{eo sum} are non-negative for fixed $j$ and sufficiently large $n$, it is necessary that there exists exactly one non-trivial linear profile $\phi^j=\phi^1$ and
\begin{align}
\psi_n(0)=T_n^1 P_n^1\phi^1+w_n^1.
\end{align}
Particularly, from \eqref{mo sum} and \eqref{eo sum} it follows
\begin{align}
\lim_{n\to\infty}\mM(T_n^1 P_n^1\phi^1)&=\mM_0,\\
\lim_{n\to\infty}\mH(T_n^1 P_n^1\phi^1)&=\mH_0,\\
\lim_{n\to\infty}\|w_n^1\|_2&=0,\label{l2 constraint w1}\\
\lim_{n\to\infty}\mH(w_n^1)&=0.\label{energy constraint w1}
\end{align}
Combining with Lemma \ref{killip visan curve} (v), \eqref{energy constraint w1} also implies
\begin{align}
\lim_{n\to\infty}\|\nabla w_n^1\|_2=0,
\end{align}
thus together with \eqref{l2 constraint w1} we deduce that
\begin{align}
\lim_{n\to\infty}\|w_n^1\|_{H^1}=0\label{h1 constraint w1}.
\end{align}
Using \eqref{orthog L2 and H1} we see that
\begin{align}
\limsup_{n\to\infty}\|T_n^1 P_n^1\phi^1\|_{H^1}<\infty.
\end{align}
From \eqref{h1 constraint w1} we infer that
\begin{align}
\limsup_{n\to\infty}\|\la\nabla \ra^{\frac{1}{2}}(\psi_n(0)-T_n^1P_n^1\phi^1)\|_{2}\lesssim\lim_{n\to\infty}\|w_n^1\|_{H^1}=0.
\end{align}
Since $\mM_0<\mM(Q)$, we obtain from \eqref{mo sum} that
\begin{align}\label{less mass}
\mM(T_n^1 P_n^1 \phi^1)\leq \frac{\mM_0+\mM(Q)}{2}<\mM(Q)
\end{align}
for sufficiently large $n$, hence Lemma \ref{large scale proxy} is applicable. If $\ld^1_\infty= \infty$, then using Lemma \ref{large scale proxy} we know that for sufficiently large $n$, there exists a global and scattering solution $v_n$ of \eqref{NLS} with
\begin{align}
v_n(0)=T_n^1P_n^1 \phi^1.
\end{align}
Finally, using Strichartz and \eqref{h1 constraint w1} we see that
\begin{align}
\lim_{n\to\infty}\|\la\nabla\ra^{\frac{1}{2}}e^{it\Delta}(\psi_n(0)-T_n^1P_n^1\phi^1)\|_{L_{t,x}^4(\R)}\lesssim\lim_{n\to\infty}\|w_n^1\|_{H^1}=0.
\end{align}
Therefore, the conditions \eqref{condition c1} to \eqref{condition a} of Lemma \ref{long time pert} are satisfied and by setting the error term $e=0$ we infer that
\begin{align}\label{contradiction4}
\limsup_{n\to\infty}\|\la\nabla\ra^{\frac{1}{2}}\psi_n\|_{{L_{t,x}^4}(\R)}<\infty,
\end{align}
which contradicts \eqref{oo1}. Hence $\ld^1_\infty=1$ and $T_n^1 P_n^1 \phi^1=T_n^1 \phi^1$. Suppose that $\lim_{n\to\infty}t_n^1=t_\infty^1$. We then define $u_c$ as the solution of the integral equation
\begin{align}
u_c(t)=e^{it\Delta}\phi^1+i\int_{-t_\infty^1}^t (|u_c|^2 u_c+|u_c|^4 u_c)(s)\,ds,
\end{align}
whose local existence near $-t_\infty^1$ is guaranteed by Lemma \ref{well posedness lemma}. \eqref{requirement1} follows already from \eqref{mo sum}, \eqref{eo sum} and the continuity of $\mD$. If \eqref{requirement2} does not hold, then by Lemma \ref{well posedness lemma} we know that $u_c$ must be a global solution. We can therefore define
\begin{align}
\tilde{v}_n:=u_c(t-t_n^1,x-x_n^1).
\end{align}
By time and space translation invariance we infer that $\tilde{v}_n$ is a solution of \eqref{NLS} and
\begin{align}
\|\la\nabla\ra^{\frac{1}{2}} \tilde{v}_n\|_{{L_{t,x}^4}(\R)}=\|\la\nabla\ra^{\frac{1}{2}} u_c\|_{{L_{t,x}^4}(\R)}.
\end{align}
By the construction of $u_c$ we also know that
\begin{align}
\lim_{n\to\infty}\|\hat{v}_n(0)-T_n^1 \phi^1\|_{H^1}=0.
\end{align}
Now as argued previously, we arrive at the contradiction \eqref{contradiction4} again. Finally, we can mimic the proof of \cite[Prop. 6.11]{Akahori2013}, words by words, to show that $u_c$ is global. We omit the details here. The proof for the first scenario is done.

We now consider the second scenario \eqref{second situation}. In this case, for each $j\in\N$ we must have
\begin{align}
\mM(T_n^jP_n^j\phi^j)&\leq \mM_0-\frac{\delta}{2}\text{ or }\nonumber\\
\mH(T_n^jP_n^j\phi^j)&\leq\mH_0-\frac{\delta}{2}
\end{align}
for sufficiently large $n$. Define
\begin{align}
\mD_1&:=\mD(\mM_0-\frac{\delta}{2},\mH_0+\vare_1(\delta)),\nonumber\\
\mD_2&:=\mD(\mM_0+\vare_2(\delta),\mH_0-\frac{\delta}{2})\nonumber\\
\end{align}
for some $\vare_1(\delta),\vare_2(\delta)>0$ such that $\mD_1,\mD_2<\mD^*$. This is possible due to Lemma \ref{killip visan curve} (iv) and the continuity of $\mD$. Thus by the inductive hypothesis \eqref{introductive hypothesis} we infer that there exist nonlinear profiles $v_n^j$ which are global solutions of \eqref{NLS} with $v_n^j(0)=T_n^jP_n^j\phi^j$ and
\begin{align}
\|\la\nabla\ra^{\frac{1}{2}}v_n^j\|_{L_{t,x}^{4}(\R)}\leq\max\{\tau(\mD_1),\tau(\mD_2)\}
\end{align}
for each $j\in\N$ and all sufficiently $n$, where $\tau$ is the quantity defined by \eqref{def of tau}. Having defined the nonlinear profiles $v_n^j$, we now define the proxy $\Psi_n^K$ by
\begin{align}
\Psi_n^K:=\sum_{j=1}^K v_n^j+e^{it\Delta}w_n^K,
\end{align}
with some sufficiently large $K$ and $n=n(K)$ to be chosen later. Since the error analysis is available for all models regardless of the signs of the nonlinearities, we are able to invoke the error analysis given in the proof of \cite[Prop. 5.2]{Cheng2020} to see that \eqref{condition c1} to \eqref{condition b} are satisfied for some sufficiently large $K$ and $n=n(K)$, where we also replace \cite[Lem. 2.3]{Cheng2020} by Lemma \ref{killip visan curve} (v). Using Lemma \ref{long time pert} we conclude the contradiction \eqref{contradiction4} again. This completes the desired proof.
\end{proof}

\section{Extinction of the minimal blow-up solution}\label{Extinction of the minimal blow-up solution}
We close in this section the proof of Theorem \ref{main theorem 1} by showing the contradiction that the minimal blow-up solution $u_c$ given by Lemma \ref{category 0 and 1} must be zero. To proceed, we still need the following lemma which can be proved as in \cite{non_radial,killip_visan_soliton} verbatim, so we omit the details of the proof here.
\begin{lemma}[\cite{non_radial,killip_visan_soliton}]\label{holmer}
Let $u_c$ be the minimal blow-up solution given by Lemma \ref{category 0 and 1}. Then
\begin{itemize}

\item[(i)] For each $\vare>0$, there exists $R=R(\vare)>0$ so that
\begin{align}
\int_{|x+x(t)|\geq R}|\nabla u_c(t)|^2+|u_c(t)|^2+|u_c|^\tas+|u_c|^\tbs\,dx\leq\vare\quad\forall\,t\in[0,\infty).
\end{align}


\item[(ii)] The center function $x(t)$ obeys the decay condition $x(t)=o(t)$ as $t\to\infty$.
\end{itemize}
\end{lemma}

\begin{proof1}
We will show the contradiction that the minimal blow-up solution $u_c$ given by Lemma \ref{category 0 and 1} is equal to zero, which will finally imply Theorem \ref{main theorem 1}. Let $\chi$ be a smooth radial cut-off function satisfying
\begin{align}
\chi(x)=\left\{
             \begin{array}{ll}
             |x|^2,&\text{if $|x|\leq 1$},\\
             0,&\text{if $|x|\geq 2$}.
             \end{array}
\right.
\end{align}
Define also the local virial action
\begin{align}
z_{R}(t):=\int R^2\chi\bg(\frac{x}{R}\bg)|u_c(t,x)|^2\,dx.
\end{align}
Direct calculation yields
\begin{align}
\pt_t z_R(t)=&\,2\,\mathrm{Im}\int R\nabla \chi\bg(\frac{x}{R}\bg)\cdot\nabla u_c(t)\bar{u}_c(t)\,dx,\label{final4}\\
\pt_{tt} z_R(t)=&\,4\int \pt^2_{jk}\chi\bg(\frac{x}{R}\bg)\pt_j u_c\pt_k\bar{u}_c-\frac{1}{R^2}\int\Delta^2\chi\bg(\frac{x}{R}\bg)|u_c|^2\nonumber\\
&\,-\int\Delta\chi\bg(\frac{x}{R}\bg)|u_c|^\tas\,dx-\frac{4}{3}\int\Delta\chi\bg(\frac{x}{R}\bg)|u_c|^\tbs\,dx.
\end{align}
Here we used the Einstein summation for the repeated indices. We then obtain that
\begin{align}\label{final 2}
\pt_{tt} z_R(t)=8\mK(u_c)+A_R(u_c(t)),
\end{align}
where
\begin{align}
A_R(u_c(t))=&\,4\int\bg(\pt^2_j\chi\bg(\frac{x}{R}\bg)-2\bg)|\pt_j u_c|^2+4\sum_{j\neq k}\int_{R\leq|x|\leq 2R}\pt^2_{jk}\chi\bg(\frac{x}{R}\bg)\pt_j u_c\pt_k\bar{u}_c\nonumber\\
&\,-\frac{1}{R^2}\int\Delta^2\chi\bg(\frac{x}{R}\bg)|u_c|^2
-\int\bg(\Delta\chi\bg(\frac{x}{R}\bg)-4\bg)|u_c|^\tas\,dx\nonumber\\
&\,-\frac{4}{3}\int\bg(\Delta\chi\bg(\frac{x}{R}\bg)-4\bg)|u_c|^\tbs\,dx.
\end{align}
We can roughly estimate $A_R(u_c(t))$ by
\begin{align}
|A_R(u_c(t))|\leq C_1\int_{|x|\geq R}|\nabla u_c(t)|^2+\frac{1}{R^2}|u_c(t)|^2+|u_c|^\tas+|u_c|^\tbs
\end{align}
for some $C_1>0$. Since $\mM(u_c)=\mM_0<\mM(Q)$, we may assume that $\mM(u_c)=(1-\delta)\mM(Q)$ for some $\delta\in(0,1)$. Using \eqref{lower bound kt} we obtain that
\begin{align}
\mK(u_c(t))&\geq\min\bg\{\delta\mH(u_c(0)),\,\bg(\bg(\frac{2}{\delta}\bg)^{\fraco}-1\bg)^{-1}\bg(m_{\mM(u_c(0))}-\mH(u_c(0))\bg)\bg\}\nonumber\\
&=: 4^{-1}\eta_1
\end{align}
for all $t\in\R$. Using Lemma \ref{holmer} (i) we know that there exists some $R_0\geq 1$ such that
\begin{align}
\int_{|x+x(t)|\geq R_0}|\nabla u_c|^2+|u_c|^2+|u_c|^\tas+|u_c|^\tbs\,dx\leq \frac{\eta}{C_1}.
\end{align}
Thus for any $R\geq R_0+\sup_{t\in[t_0,t_1]}|x(t)|$ with some to be determined $t_0,t_1\in[0,\infty)$, we obtain that
\begin{align}\label{final3}
\pt_{tt} z_R(t)\geq \eta_1
\end{align}
for all $t\in[t_0,t_1]$. By Lemma \ref{holmer} (ii), for some to be determined small $\eta_2$ we can choose sufficiently large $t_0$ such that $|x(t)|\leq\eta_2 t$ for all $t\geq t_0$. Now set $R=R_0+\eta_2 t_1$. Integrating \eqref{final3} over $[t_0,t_1]$ yields
\begin{align}\label{12}
\pt_t z_R(t_1)-\pt_t z_R(t_0)\geq \eta_1 (t_1-t_0).
\end{align}
Using \eqref{final4}, Cauchy-Schwarz and Lemma \ref{killip visan curve} (v) we have
\begin{align}\label{13}
|\pt_t z_{R}(t)|\leq C_2 \mD^*R= C_2 \mD^*(R_0+\eta_2 t_1)
\end{align}
for some $C_2=C_2(\mD^*)>0$. \eqref{12} and \eqref{13} give us
\begin{align}
2C_2 \mD^*(R_0+\eta_2 t_1)\geq\eta_1 (t_1-t_0).
\end{align}
Setting $\eta_2=(4C_2\mD^*)^{-1}$ and then sending $t_1$ to infinity we obtain a contradiction unless $\eta_1=0$, which implies $\mH_0=\mH(u_c)=0$. From Lemma \ref{killip visan curve} (v) we conclude that $\nabla u_c=0$, which implies $u_c=0$. This completes the proof.
\end{proof1}

\section{A solely mass-determining scattering threshold for the problem \eqref{aka nls}}\label{outlook section}
In this section we continue our discussion on the NLS \eqref{aka nls}. Our aim is to impose a solely mass-determining scattering threshold for \eqref{aka nls} which is similar to the one given in Theorem \ref{main theorem 1}. Our starting point is the following result given in \cite{SoaveCritical,wei2021normalized}:
\begin{theorem}[\cite{SoaveCritical,wei2021normalized}]
Let $d\geq 3$. Then
\begin{itemize}
\item[(i)]\textbf{Existence of ground state}: For any $c\in(0,\infty)$ the variational problem \eqref{soave problem} corresponding to \eqref{aka nls} has a minimizer $S_c$ with $\mH(S_c)=m_c\in(0,d^{-1}\csob^{\frac{d}{2}})$. Moreover, $S_c$ is a solution of \eqref{standing wave eq} with some $\omega>0$. In addition, $S_c$ can be chosen to be positive and radially symmetric.
\item[(ii)]\textbf{Blow-up criterion}: Assume that $u_0\in H^1(\R^d)$ satisfies the conditions $\mH(u_0)< m_{\mM(u_0)}$ and $\mK(u_0)<0$. Assume also that $|x|u_0\in L^2(\R^d)$. Then the solution $u$ of \eqref{aka nls} with $u(0)=u_0$ blows-up in finite time.
\end{itemize}
\end{theorem}

The scattering result can now be formulated as follows:
\begin{theorem}\label{3d scattering threshold}
Define the set
\begin{align}
\mathcal{B}&:=\{u\in H^1(\R^d):\mH(u)<m_{\mM(u)},\mK(u)> 0\}
\end{align}
and assume that $u_0\in \mathcal{B}$. Additionally we assume that $u_0$ is radially symmetric in the case $d=3$. Then the solution $u$ of \eqref{aka nls} with $u(0)=u_0$ is global and scatters in time.
\end{theorem}
The proof is a straightforward modification and combination of the variational arguments given in Section \ref{Variational problems and estimates} and the nonlinear estimates in \cite{Akahori2013}, we therefore omit the details here.

\subsubsection*{Acknowledgments}
The author acknowledges the funding by Deutsche Forschungsgemeinschaft (DFG) through the Priority Programme SPP-1886 (No. NE 21382-1). The author also thanks the anonymous referee sincerely for her/his thorough reading of the manuscript and for many important corrections.

\begin{appendix}
\section{Precise values of $m_0$ and $m_Q$}\label{m0 mQ}
\begin{proposition}
We have $m_0=\infty$ and $m_Q=0$.
\end{proposition}

\begin{proof}
From Theorem \ref{soave} we know that for any $c\in(0,\mM(Q))$, $m_c$ has a minimizer $S_c$. Therefore
\begin{align}
\mH(S_c)=m_c,\quad\mK(S_c)=0,\quad\mM(S_c)=c.
\end{align}
Using $\mK(S_c)=0$ and Gagliardo-Nirenberg we infer that
\begin{align}
\|\nabla S_c\|_2^2=\frac{1}{2}\|S_c\|_4^4+\frac{2}{3}\|S_c\|_6^6\leq c \mM(Q)^{-1}\|\nabla S_c\|_2^2 +\frac{2 c\,\gnb^{-1}}{3}\|\nabla S_c\|_2^4,
\end{align}
which implies
\begin{align}
\|\nabla S_c\|_2^2\geq \frac{3\gnb}{2}\bg(\frac{1}{c}-\frac{1}{\mM(Q)}\bg).
\end{align}
Therefore $\|\nabla S_c\|_2^2\to\infty$ as $c\to 0$. Now we obtain
\begin{align}
\mH(S_c)&=\mH(S_c)-\frac{1}{4}\mK(S_c)\nonumber\\
&=\frac{1}{4}(\|\nabla S_c\|_2^2-\frac{1}{2}\|S_c\|_4^4)\nonumber\\
&\geq\frac{1}{4}(1-c\mM(Q)^{-1})\|\nabla S_c\|_2^2\to \infty
\end{align}
as $c\to 0$, which implies $m_0=\infty$. Next we show $m_Q=0$. Let $(u_n)_n$ be a minimizing sequence for \eqref{def gn l2 crit}. By rescaling we may assume that $\mM(u_n)=\delta^2 \mM(Q)$ for $\delta\in(1/2,1)$ which will be sended to one later, and $\|u_n\|_4\equiv 1$. Then combining with \eqref{GN-L1} we obtain that $\|\nabla u_n\|_2^2=(2\delta^2)^{-1}+o_n(1)$. We then conclude that
\begin{align}
\mK(T_\ld u_n)= \frac{\ld^2}{2}\bg(\frac{1}{\delta^2}-1+o_n(1)\bg)-\frac{2\ld^4}{3}\|u_n\|_6^6.
\end{align}
By setting
\begin{align}
\ld_{n,\delta}=\bg(\frac{3}{4\|u_n\|_6^6}\bg(\frac{1}{\delta^2}-1+o_n(1)\bg)\bg)^{\frac{1}{2}}
\end{align}
we see that $\mK(T_{\ld_{n,\delta}}u_n)=0$. By H\"older we obtain that
\begin{align}
\|u_n\|_6^6\geq \mM(u_n)^{-1}\|u_n\|_4^8=\delta^{-2}\mM(Q)^{-1}.
\end{align}
We now choose $N=N(\delta)\in\N$ such that $|o_n(1)|\leq \delta^{-2}-1$ for all $n>N$. Summing up and using the definition of $m_c$ we finally conclude that
\begin{align}
m_{\delta^2 \mM(Q)}&\leq \sup_{n>N}\mH(T_{\ld_{n,\delta}}u_n)=\sup_{n>N}\bg(\mH(T_{\ld_{n,\delta}}u_n)-\frac{1}{2}\mK(T_{\ld_{n,\delta}}u_n)\bg)\nonumber\\
&=\sup_{n>N}\frac{1}{6}\|T_{\ld_{n,\delta}}u_n\|_6^6=\sup_{n>N}
\frac{\ld_{n,\delta}^4}{6}\|u_n\|_6^6\nonumber\\
&\leq \frac{9}{24}\mM(Q)\delta^{-2}(1-\delta^2)^2\to 0
\end{align}
as $\delta\to 1$. This proves $m_Q=0$.
\end{proof}
\end{appendix}

\addcontentsline{toc}{section}{References}

\end{document}